\documentclass[pdflatex,sn-mathphys-num]{sn-jnl}

\usepackage{graphicx}%
\usepackage{multirow}%
\usepackage{amsmath,amssymb,amsfonts}%
\usepackage{amsthm}%
\usepackage{mathrsfs}%
\usepackage[title]{appendix}%
\usepackage{xcolor}%
\usepackage{textcomp}%
\usepackage{manyfoot}%
\usepackage{booktabs}%
\usepackage{algorithm}%
\usepackage{algorithmicx}%
\usepackage{algpseudocode}%
\usepackage{listings}%
\usepackage{mathtools}
\usepackage{subcaption}
\usepackage{hyperref}

\newcommand{\dt}{{\Delta t}}

\newcommand{\s}{\sigma}
\renewcommand{\a}{\alpha}
\renewcommand{\b}{\beta}
\newcommand{\ddt}{\partial_t}
\newcommand{\ddx}{\partial_x}
\newcommand{\ddy}{\partial_y}

\newcommand{\ket}[1]{| #1 \rangle}

\newcommand{\lenv}[1]{\textbf{L}_{\a_{#1-1}}^{(#1)}}

\newcommand{\renv}[1]{\textbf{R}_{\a_{#1}}^{(#1)}}

\newcommand{\isel}[1]{\mathcal{I}^{(#1)}_{\a_{#1-1}}}
\newcommand{\jsel}[1]{\mathcal{J}^{(#1)}_{\a_{#1}}}
\newcommand{\bisel}[1]{\mathcal{I}^{(#1)}_{\a'_{#1-1}}}
\newcommand{\bjsel}[1]{\mathcal{J}^{(#1)}_{\a'_{#1}}}

\theoremstyle{thmstyleone}%
\theoremstyle{thmstyletwo}%

\theoremstyle{thmstylethree}%

\begin{document}

\title[Article Title]{Time integration of quantized tensor trains using the interpolative dynamical low-rank approximation}


\author*[1]{\fnm{Erika} \sur{Ye}}\email{erikaye@lbl.gov}

\author[1]{\fnm{Chao} \sur{Yang}}\email{cyang@lbl.gov}

\affil*[1]{\orgdiv{Applied Mathematics and Computational Research Division}, \orgname{Lawrence Berkeley National Laboratory}, \orgaddress{\street{1 Cyclotron Rd}, \city{Berkeley}, \postcode{94720}, \state{CA}, \country{USA}}}


\abstract{Quantized tensor trains (QTTs) are a low-rank and multiscale framework that allows for efficient approximation and manipulation of multi-dimensional, high resolution data. One area of active research is their use in numerical simulation of hyperbolic systems such as the Navier-Stokes equations and the Vlasov equations. One popular time integration scheme is the dynamical low-rank approximation (DLRA), in which the time integration is constrained to a low-rank manifold. However, until recently, DLRA has typically used orthogonal projectors to project the original dynamical system into a reduced space, which is only well-suited for linear systems. DLRA has also mostly been investigated in the context of non-quantized tensor trains. This work investigates interpolative DLRA schemes in which the low-rank manifold is constructed from aptly chosen interpolation points and interpolating polynomials, in the context of QTTs. Through various examples, its performance is compared to its orthogonal counterpart. This work demonstrates how interpolative DLRA is suitable for nonlinear systems and time integrators requiring nonlinear element-wise operations, such as upwind time integration schemes.}




\maketitle

\section{Introduction}\label{sec1}

Tensor trains (TTs) are a powerful low-rank numerical framework that has been used to solve eigenvalue problems (e.g., in quantum chemistry applications) and initial value problems (e.g., classical time-dependent PDEs) with significantly reduced cost. 
Traditionally, the tensor train is interpreted as a weighted sum of an incomplete set of orthonormal basis functions. These basis functions are obtained naturally from the SVD and QR decompositions used to construct the TT. This interpretation will be referred to as the orthonormal construction.
However, many recent works have begun utilizing an interpolative construction, in which the state is interpreted as a sum of interpolating polynomials with some choice in interpolation nodes. The underlying decomposition is the skeleton or CUR decomposition.
Unlike in the orthonormal construction, the interpolative construction provides element-wise access and thus enables efficient nonlinear evaluations.

The focus of this paper is on using quantized tensor trains (QTTs) for time integration of classical partial differential equations. The QTT is a low-rank ansatz designed to take advantage of any multi-scale structure,  potentially allowing one to solve PDEs with reduced cost. For a $K$-dimensional problem with $N$ grid points along each dimension, the computational cost is reduced from $\mathcal{O}(N^K)$ to $\mathcal{O}(r^3 K \log N)$, where $r$ is the rank of the QTT. Though $r$ will depend on the problem at hand, it should be independent with respect to $N$, if the dynamics is adequately resolved.  The various algorithms for QTT time integration can be sorted into two broad categories: global update methods and methods based on the dynamical low-rank approximation. 
Both the orthonormal and interpolative TT constructions can be used in either class of time integration algorithms. Table \ref{tab:summary} provides a high-level summary of different time evolution schemes. 

In global update schemes, one can use any existing time integration method to update the TT as a single entity, with all internal steps replaced by their TT analogs. 
Simple operations like vector addition and matrix-vector multiplication are performed using the step-and-truncate (SAT) procedure, where an exact TT representation of the result is obtained and then a rank-truncation procedure is performed \cite{schollwock_density-matrix_2011}. One can also solve linear equations using an alternating least-squares optimization \cite{Oseledets2012solution, holtz_alternating_2012} or an enriched version such as the alternating minimal energy algorithm \cite{dolgov_tensor_2014, Dolgov2014alternating}. These methods are most commonly used with TTs in the orthonormal construction. However, they have been adapted to work with the interpolative construction to tackle problems requiring nonlinear function evaluations \cite{dolgov_polynomial_2015}. Ref.~\cite{appelo_lraa_2025} also recently introduced a cross algorithm to use Anderson acceleration for solving nonlinear matrix equations. 
Similar iterative algorithms can be used more generally for efficiently evaluating nonlinear operations on TTs. Example situations include performing element-wise multiplication with a irregular function \cite{peddinti_complete_2023}, evaluating nonlinear contributions in the PDE \cite{dolgov_polynomial_2015, gharemani_cross_2024}, or applying filtering in higher-order time integration methods \cite{danis_weno_2024}. 
Recent works even demonstrate performing semi-Lagrangian time integration by using an adaptive cross approximation procedure to identify the reduced set of points at which one needs to compute the characteristic curves to obtain the solution at the next time step \cite{zheng_SLAR_2024, zheng_SLAR_2025, christlieb_sampling_2025}.

In the dynamical low-rank approximation (DLRA), the dynamics of each tensor core is approximated by projecting the original equation of motion onto the manifold defined by the remaining tensor cores of the tensor train \cite{koch_dynamical_2007}. 
In practice, a projector splitting (PSI) split-step implementation is used and each tensor core is updated one at a time \cite{lubich_projector-splitting_2014}, as the original formulation is nonlinear and prone to numerical instabilities. This often results in reduced costs in practice, though the computational complexity is formally the same as in the SAT procedure. Additionally, its straightforward to use a wide range of time integrators to solve the reduced problem at each tensor core.
Further modifications in the details of the implementation and the ordering in which they are applied have been investigated. For example, in the Galerkin basis update (BUG) \cite{ceruti_robust_2024}  and the Galerkin alternating projection (GAP) methods \cite{ceruti_galerkin_2025}, the backwards time steps in the PSI algorithm are replaced by projections of the original TT onto the updated manifold. This improves stability in the presence of dissipative terms. The algorithm can also easily be modified to allow for rank adaptivity \cite{ceruti_rank-adaptive_2022}. Ref.~\cite{einkemmer_review_2024} provides a thorough review of these methods.
With the orthonormal construction, DLRA can be used to solve various types of differential equations \cite{hochbruck_rank-adaptive_2023}, and there already exist many numerical demonstrations for a variety of applications, including the kinetic equation \cite{einkemmer_low-rank_2020, ceruti_robust_2024}, the gyrokinetic equations \cite{einkemmer_accelerating_2023}, and radiation transport \cite{einkemmer_asymptotic-preserving_2024}. 
Though many of the above references consider matrix data sets (i.e., tensor trains of length two), they can be easily extended to the tensor train case \cite{lubich_dynamical_2013}. 
Additionally, the projector splitting time integrator has also been used in various quantum chemistry and condensed matter applications, though they are instead referred to as time-dependent density matrix renormalization group (TD-DMRG) \cite{feiguin2005time} and the time-dependent variational principle (TDVP) \cite{Haegeman_TDVP}. Ref.~\cite{Paeckal_2019_time} provides a detailed overview of these methods in this context.

Studies using the interpolative construction are more limited. 
Ref. \cite{dektor_collocation_2024} introduces this method for tensor trains and provides a few numerical examples. 
Ref.~\cite{dektor_interpolatory_2024} studies the Boltzmann-BGK equation. Because of dissipative nonlinear term, the original method is modified to replace the backwards time integration with a basis update step, arriving at the DLR-basis update and collate (DLR-BUC) scheme. However, this method was only discussed in the context of matrix decompositions (ie., tensor trains of length two). 
Ref.~\cite{hossein_2025_cross} developed an implicit DLR-based solver for a nonlinear PDEs with parametric randomness, though it also is presented for the matrix case.  

The goal of this work is to investigate the practicality of interpolative DLR-type time integration in the context of \textit{quantized} TTs, which has not yet been done.
The mechanics of the time integrator can be easily extrapolated from the standard DLRA method for QTTs, as well as the interpolative DLRA method for functional TTs.
However, conceptually, the interpolative construction may seem odd, particularly in the quantized setting. For example, in contrast to most functional TT decompositions, operators such as the spatial derivative are not separable across the quantized dimensions. Additionally, while the CUR decomposition can yield a near-optimal matrix decomposition \cite{savostyanov_quasioptimality_2014}, its compatibility with DLRA framework, and the impact of the error arising from projecting the dynamics onto the interpolative TT manifold, is unclear. 

The paper is structured as follows. The first section provides the necessary background information on quantization, tensor trains and the interpolative construction, and DLRA-type time integrators. Analogies between the interpolative and the orthogonal construction are highlighted, so that the extension from classical DLR schemes to interpolative DLR schemes should be straightforward. This is followed by a deeper introduction into two ways interpolative DLRA can be performed--by using the interpolative TT construction, and by using oblique projections to convert the TT between an interpolative and orthonormal representation.
This is then followed by some simple numerical examples with QTTs. The first example is the time evolution of invsicid Burger's equation with upwinding solved using interpolative DLRA. 
The second example is an electromagnetic pulse in a dielectric cavity, in which the various DLRA-type integrators and the global-update (step-and-truncate) time integration are compared. The last example is a simple advection problem, in which higher order time integration is demonstrated. Furthermore, the results are compared between different types of DLRA integrators and the benefits of using QTTs for calculations with increased resolution are clearly shown. Finally, the paper concludes with a summary of the above results and a discussion for future research directions. 

\begin{table}[h]
    \caption{List of examples for each type of time integrator in the tensor train format.}
    \label{tab:summary}
    \begin{tabular} {|p{1.6cm}|p{5cm}|p{5cm}|}
        \hline  
        & \textbf{Global update schemes} & \textbf{Local update schemes} \\
        \hline  
        \vspace{0.5cm}
        Orthogonal & 
        \vspace{-0.2cm}
        \begin{itemize}
            \item Step-and-truncate
            \vspace{-0.1cm}
            \begin{itemize}
                \setlength\itemsep{0.0em}
                \item spin models \cite{vidal_efficient_2004}
                \item Maxwell's equations \cite{manzini_tensor_2023}
                \item Vlasov equation \cite{guo_conservative_2024, ye_quantum-inspired_2022}
            \end{itemize}
            \item DMRG-based optimzation or solving
            \vspace{-0.1cm}
            \begin{itemize}
                \setlength\itemsep{0.0em}
                \item Navier-Stokes \cite{gourianov2022exploiting}) 
                \item Chemical master equation \cite{kazeev_direct_2014}
            \end{itemize}
        \end{itemize}
        &
        \vspace{-0.2cm}
        \begin{itemize}
            \item DLRA-PS or TDVP
            \vspace{-0.1cm}
            \begin{itemize}
                \setlength\itemsep{0.0em}
                \item spin models \cite{Haegeman_TDVP}
                \item Vlasov  \cite{einkemmer_low-rank_2020, ye_quantized_2023} 
                \item Gyrokinetic equation \cite{einkemmer_accelerating_2023}
            \end{itemize}
            \item DLRA-AP or TD-DMRG, DLRA-BUG 
            \vspace{-0.1cm}
            \begin{itemize}
                \setlength\itemsep{0.0em}
                \item Quantum chemistry \cite{feiguin2005time, Ronca_TDDMRG} 
                \item Heat equation, Vlasov-Poisson \cite{ceruti_robust_2024}
                \item Gyrokinetic equation \cite{einkemmer_accelerating_2023}
            \end{itemize}
        \end{itemize}
        \\
        \hline
        \vspace{0.3cm}
        Interpolative  & 
        \vspace{-0.2cm}
        \begin{itemize}
            \item nonlinear PDEs \cite{gharemani_cross_2024}
            \item upwinding and flux limiters \cite{danis_weno_2024}
            \item non-uniform masks \cite{peddinti_complete_2023}
            \item Semi-Lagrangian time evolution \cite{zheng_SLAR_2024, zheng_SLAR_2025, christlieb_sampling_2025}
            \item Anderson acceleration \cite{appelo_lraa_2025}
        \end{itemize}
        &
        \vspace{-0.2cm}
        \begin{itemize}
            \item DLR with collocation \cite{dektor_collocation_2024}
            \item DLR-BUC: Boltzmann-BGK model \cite{dektor_interpolatory_2024}
            \item Parametric PDEs
            \cite{hossein_2025_cross, budzinsky_low-rank_2025}
            \item DLR with Tucker, random index selection \cite{carrel_interpolatory_2025}
            \item this paper 
        \end{itemize}
        \\
        \hline  
    \end{tabular}
\end{table}



\section{Background}

\subsection{Quantization}
Consider some arbitrary 1-dimensional function $f$ expanded with respect to some basis $\{ {\phi}_n \} $, 
\begin{align}
    \ket{f(x)} = \sum_{n} f_{n} \ket{ \phi_n(x) } 
\end{align}
where $f_n = \langle \phi_n | f \rangle$. 
Quantization is a folding \cite{Khoromskij2011} or mapping procedure \cite{Lindsey2023multiscale} in which the one-dimensional function is converted into a multi-dimensional representation. 

Suppose the basis is of size $N=2^L$.
In the quantized representation, one defines a one-to-one mapping of each integer $n$ to hyperindex $(\sigma_1, \sigma_2, ... \s_L)$ where $\sigma_p \in \{0, 1 \}$.
For example, let $f$ be evaluated in real space $x\in[0,1)$ at $N$ uniform grid points. In this case, $\phi_n(x) = \delta(x-x_n)$ for $n \in \{0,1,...N-1\}$, where $x_n = n/N$. Then one can define a mapping $n \leftrightarrow (\sigma_1, \s_2, \hdots \s_L)$ with
\begin{align}
    x_n = 2^{-L} \, n = \sum_{k=1}^\infty  2^{-k} \sigma_k = 0.\sigma_1 \sigma_2 \sigma_3 \hdots \sigma_L
\end{align}
where the expression on the right-hand side denotes a binary decimal expansion.
%
%
From this mapping, one then obtains the {quantized tensor representation} of the coefficients $f_n$,
\begin{align}
    \ket{f(x)} = \sum_{n=1}^N f_n \ket{\phi_n(x)} 
    \cong 
    \sum_{\s_1} \sum_{\s_2} \hdots \sum_{\s_L}  T(\s_1,\s_2,\hdots,\s_L) \,\ket{\phi_{\s_1,\s_2,\hdots, \s_L} (x)}.
\end{align}
where $T$ is an $L$-dimensional tensor of size $2$ along each dimension.
%
With this choice of basis, the $\{ \sigma_i\}$ indices correspond to grid scale, and the grids at different grid scales are separable. 
For example, 
the grid can be partitioned at depth $m$ into a coarse grid and fine grid,
\begin{align}
    x_{\leq m} := \sum_{k=1}^m 2^{-k} \sigma_k, 
    \qquad
    x_{>m} := \sum_{k=m+1}^L 2^{-k} \sigma_k = 2^{-m} \sum_{k=1}^{L-m} 2^{-k} \sigma_k.
\end{align}
%
%
While this property is not necessary for using the methodology that will be discussed in the remainder of the paper, its multi-scale nature motivates the use of QTTs for solving multi-scale partial differential equations (PDEs) that are expected to exhibit some amount of scale separation, in which the behavior at each grid scale only strongly depends on the behavior at adjacent grid scales. This work does not study multi-scale systems, but instead uses simpler test problems for proof of concept.

\subsection{Tensor train format}

\subsubsection{Vectors}

Consider the $L$-dimensional tensor $T(\sigma_1, \sigma_2, \hdots, \sigma_K) \in \mathbb{C}^{d_1\times d_2 \times ... \times d_L}$. 
The tensor can be decomposed into the tensor train (TT) format, 
\begin{align}
    T(\sigma_1, \hdots \sigma_L) = 
    \sum_{\alpha_1 = 1}^{r_1} \hdots
    \sum_{\alpha_{L-1} = 1}^{r_{L-1}} 
    \left(
    M^{(1)}_{\alpha_1} (\sigma_1) \, M^{(2)}_{\alpha_1, \alpha_2} (\sigma_2) \hdots M^{(L-1)}_{\alpha_{L-2}, \alpha_{L-1}} (\sigma_{L-1}) \, M^{(L)}_{\alpha_{L-1}}(\sigma_L)
    \right)
\end{align}
consisting of $L$ tensor cores $M^{(i)}$ of size $r_{i-1} \times d_i \times r_{i}$ (where $r_0=r_L = 1$ and thus can be dropped). The rank of the tensor train is $r = \max(r_1, \hdots r_{L-1})$. The tensor notation is chosen such that the dimension written in parentheses will correspond to a physical dimension (which have physical meaning and still remain after contraction of all tensors in the tensor train) while the subscripts will correspond to virtual dimensions (which arise from the tensor decomposition and are summed over during tensor contraction).

Due to gauge invariance, there are infinitely many ways one can choose each tensor core $M^{(i)}$ for a given tensor $T$. However, for computational efficiency and accuracy, the TT-vector is often kept in ``canonical form'', where one tensor core is the center of orthogonality and the tensors to its left/right are in ``left/right canonical form'': 
\begin{align}
    T(\sigma_1, \hdots \sigma_L) = 
    \sum_{\alpha_1 = 1}^{r_1} \hdots
    \sum_{\alpha_{L-1} = 1}^{r_{L-1}} 
    & \left(
    U^{(1)}_{\alpha_1} (\sigma_1) \, \hdots U^ {(i-1)}_{\alpha_{i-2}, \alpha_{i-1}} (\sigma_{i-1}) \right)
    \, 
    \times
    \nonumber
    \\
    & M^{(i)}_{\alpha_{i-1}, \alpha_i} (\sigma_i) 
    \left( V^{(i+1)}_{\alpha_{i}, \alpha_{i+1}} (\sigma_{i+1}) \hdots \, V^{(L)}_{\alpha_{L-1}}(\sigma_L)
    \right)
\end{align}
where $U^{(i)}$ and $V^{(i)}$ are the tensor cores in left and right canonical form, respectively. Tensor core $M^{(i)}$ is an arbitrary tensor and is the center of orthogonality.
For notational simplicity, the left and right canonical tensors will be collectively represented as a single tensor, 
\begin{align}
    \textbf{L}^{(i)}_{\alpha_{i}}(\s_{1:i-1}) &= \sum_{\a_1, \hdots, \a_{i-1}} \, U_{\a_1}^{(1)}(\s_1) \, U_{\a_1,\a_2}^{(2)}(\s_2) \, \hdots \, U_{\a_{i-1},\a_{i}}^{(i)} (\s_{i-1})
    \\
    \textbf{R}^{(i)}_{\alpha_{i+1}}(\s_{i+1:L}) &= \sum_{\a_{i+2}, \hdots, \a_{K-1}} V_{\a_{i+1},\a_{i+2}}^{(i+1)}(\s_{i+1}) \, V_{\a_{i+2},\a_{i+3}}^{(i+2)}(\s_{i+2}) \, \hdots \, V_{\a_{L-1}}^{(L)} (\s_{L})
\end{align}
where we use $\sigma_{i:j}$ as shorthand for the collection of indices $(\s_i, \hdots, \s_j)$. 
The tensor train can now be written as 
\begin{align}
    T(\s_{1:L}) = \sum_{\a_{i-1},\a_{i}} \textbf{L}^{(i)}_{\a_{i-1}} (\s_{1:i-1}) \, M^{(i)}_{\a_{i-1},\a_{i}} (\s_i) \, \textbf{R}^{(i)}_{\a_{i}} (\s_{i+1:L})
\end{align}
which we will often write as $\textbf{L} M \textbf{R}$ for short. 
$\textbf{L}$ and $\textbf{R}$ also define the manifold on which the tensor train lives, acting as subspaces of sizes $r_i$ and $r_{i+1}$ on which the original tensor is projected onto. Following this interpretation, the tensor can be written as 
\begin{align}
    T(\s_{1:L}) = \sum_{\a_{i-1},\a_{i}} M^{(i)}_{\a_{i-1},\a_{i}, \s_i}
    \ket{ \textbf{L}^{(i)}_{\a_{i-1}} \, \phi^{(i)}_{\s_{i}}, \, \textbf{R}^{(i)}_{\a_{i}} }
    \label{eq:lmr}
\end{align}
where the ket represents the canonical vectors constructing the tensor train manifold.

In the standard tensor train construction, the canonical form requires that the  vectors $|\textbf{L}^{(i)}_{\a_{i-1}},\phi^{(i)}_{\s_i}, \textbf{R}^{(i)}_{\a_{i}} \rangle$ are a set of orthonormal basis functions, and the center of orthogonality $M^{(i)}$ corresponds to the overlap of the original tensor with each basis function. 
This is satisfied if the canonical tensor cores $U^{(i)}$ and $V^{(i)}$ satisfy the orthonormality conditions,
\begin{align*}
    & \sum_{\alpha_{i-1}, \, \sigma_i} U^{(i)*}_{\alpha_{i-1}, \alpha_i'} (\sigma_i) U^{(i)}_{\alpha_{i-1}, \alpha_i} (\sigma_i) = \delta_{\alpha_i', \alpha_i} 
    & \hspace{2cm} \text{(left canonical form)}
    \\
    & \sum_{\alpha_{i}, \, \sigma_i} \, V^{(i)*}_{\alpha_{i-1}', \alpha_i} (\sigma_i) V^{(i)}_{\alpha_{i-1}, \alpha_i} (\sigma_i) = \delta_{\alpha_{i-1}', \alpha_{i-1}} 
    & \hspace{2cm} \text{(right canonical form)}.
\end{align*}
The matrices $U$ and $V$ are typically obtained from the singular vectors of the SVD, since the error of the decomposition is given by the singular values corresponding to the truncated singular vectors \cite{Oseledets2011}. In the case of infinite rank, the basis vectors span the full computational space. 

For TTs in the interpolative construction, a different definition of canonical form is used. Here, the basis vectors are interpolating functions, and the center of orthogonality contains the values of the original tensor evaluated at the interpolation points \cite{Oseledets2010cross, savostyanov_fast_2011, savostyanov_quasioptimality_2014}. 
In this case, the tensor cores are in canonical form if 
\begin{align}
    U^{(i)}_{(\alpha_{i-1},\sigma_i),\alpha_i}[I^{(i)}_{\alpha_i'},:] &=  
    \delta_{\alpha_i', \alpha_i} & \text{(left canonical form)}
    \label{eq:canon_cx}
    \\
    V^{(i)}_{\alpha_{i-1},(\alpha_i,\sigma_i)} [:, J^{(i)}_{\alpha_{i-1}'}] &= 
    \delta_{\alpha_{i'}, \alpha_{i-1}} & \text{(right canonical form)} 
    \label{eq:canon_xr}
\end{align}
where $I^{(i)}$ and $J^{(i)}$ are sets of row and column indices specific to the $i^\text{th}$ tensor core.
These tensors are typically obtained using the CUR decomposition, $M=C \hat{M}^{-1} R$, where $C=M[:,J]$ are the sub-selected columns, $R=M[I,:]$ is the sub-selected rows, and $\hat{M}$ is the submatrix $M[I,J]$. The selection indices $I$ and $J$ are determined via some column or row sub-selection algorithm. Here, the maxvol algorithm \cite{goreinov_how_2010} or the q-DEIM algorithm were used \cite{drmac_new_2016}. 
If $\hat{M}$ is non-singular, then with $U \leftarrow C\hat{M}^{-1}$ and $V \leftarrow \hat{M}^{-1} R$, the above conditions are satisfied.  


However, if $\hat{M}$ is singular, which would occur if the matrix is oversampled and the number of row or column indices is greater than the rank of the $M$, then the above conditions are not guaranteed. Instead, one can consider a set of relaxed constraints, 
\begin{align}
    U^{(i)}_{(\alpha_{i-1},\sigma_i),\alpha_i}[I^{(i)}_{\alpha_i'},:] R^{(i)} &=
    R^{(i)} & \text{(left canonical)} 
    \label{eq:canon_xr_weak}
    \\
    C^{(i)} V^{(i)}_{\alpha_{i-1},(\alpha_i,\sigma_i)} [:, J^{(i)}_{\alpha_{i-1}'}] &= 
    C^{(i)}  & \text{(right canonical)} 
    \label{eq:canon_cx_weak}
\end{align}
which is automatically satisfied by the CUR decomposition.
However, if one updates the rows $R$ in Eq.~\eqref{eq:canon_xr_weak} or columns $C$ in Eq.~\eqref{eq:canon_cx_weak} while holding the interpolating matrix fixed, as is commonly done in alternating DMRG-style algorithms, then the above conditions are generally no longer necessarily satisfied. 
The one exception is if the rows or columns are only scaled by some arbitrary scalar constant.

One can again show that the collective states $\textbf{L}$ and $\textbf{R}$ are in left and right canonical form through an iterative index selection procedure starting from the ends of the tensor train. For example, for a block built from two left canonical sites,
\begin{align*}
     \sum_{\alpha_1} \left( U^{(1)}_{\sigma_1,\alpha_1} U^{(2)}_{\alpha_1,\sigma_2,\alpha_2} \right) [\mathcal{I}^{(2)}_{\alpha_2'}] \, R_{\alpha_2, \hdots} &= 
     \sum_{\alpha_1} U^{(1)}_{\sigma_1,\alpha_1} [I^{(1)}_{\alpha_1'},:] \, \left( U^{(2)}_{(\alpha_1,\sigma_2),\alpha_2} [I^{(2)}_{\alpha_2'},:] \, R_{\alpha_2, \hdots} \right)
     \nonumber \\
     &
     = U^{(2)}_{(\alpha'_1,\sigma_2),\alpha_2}  [{I^{(2)}_{\alpha_2'}}, :] \, R_{\alpha_2, \hdots} = R_{\alpha_2, \hdots}
\end{align*}
The selection indices for $\textbf{L}$ and $\textbf{R}$ are the collection of selection indices at each tensor core, $\mathcal{I}^{(i)}=\{ I^{(1)},I^{(2)},\hdots, I^{(i-1)} \}$ and $\mathcal{J}^{(i)} = \{ J^{(i+1)}, J^{(i+2)}, \hdots, J^{(L)} \}$. The indices are \emph{nested} \cite{savostyanov_fast_2011, savostyanov_quasioptimality_2014}, since the $I^{(i)}$ selection indices only consider the rows that were selected by the $I^{(i-1)}$ selection indices. Likewise, the $J^{(i)}$ selection indices only consider the columns selected by the $J^{(i+1)}$ indices.



To convert an arbitrary TT into the canonical interpolative form \cite{savostyanov_fast_2011}, one can follow a procedure analogous to that of the orthogonal case \cite{Oseledets2011}, in which one starts from two ends of the tensor train and performs the CUR decomposition repeatedly to obtain the left and right canonical tensor cores.

\subsubsection{Linear operators}
Linear operators acting on the $L$-dimensional space can also be written in tensor train form,
\begin{align}
    A(\s'_{1:L}, \, \s_{1:L}) = 
    \sum_{\b_i=1}^{r_1} \hdots \sum_{\b_{K-1}=1}^{r_{K-1}} & \biggl( M_{\b_1}^{(1)}(\s'_1,\s_1) M_{\b_1, \b_2}^{(2)}(\s'_2,\s_2) \hdots
    \nonumber 
    \\
    & \,
    M^{(L-1)}_{\b_{L-2},\b_{L-1}} (\s'_{L-1},\s_{L-1}) M^{(L)}_{\b_{L-1}} (\s'_L,\s_L) \biggr)
\end{align}
where the hyperindices $\s_{1:L}$ and $\s_{1:L}'$ correspond to the input and output and the operator, respectively. 


Naively, matrix-vector multiplication involves contracting the $\{\s_i\}$ indices between the corresponding tensor cores of the QTT-operator and QTT-vector, such that the output retains the tensor train form \cite{dolgov_tensor_2014, schollwock_density-matrix_2011}. However, the resulting tensor train will have an enlarged rank, and one must use a rank truncation procedure to reduce the rank down to some acceptable value \cite{Oseledets2011}. However, the naive algorithm of first using QR decompositions to ensure the TT is in canonical form and then using SVDs to perform the low-rank approximation is quite slow. The zip-up algorithm \cite{stoudenmire_minimally_2010} offers noticeable speed-up. Recently, Ref.~\cite{camano_successive_2025} introduced a randomized algorithm that further reduces the cost of the this operation, though it was not used in this work. 

In the QTT format, structured matrices often have relatively low-rank representations. For example, circulant and Toeplitz matrices will have QTT of rank $2p$, where $p$ is the QTT rank of the vector constructing the matrices \cite{kazeev_multilevel_2013}. It follows that Jordan matrices are of rank two, while tridiagonal Toeplitz matrices are of rank 3 \cite{kazeev_low-rank_2012, dolgov_tensor_2014}. This implies that differential operators on a uniform grid will also be low-rank in the QTT format.

\subsection{Projector splitting time integrator}

In the dynamical low-rank scheme, the time derivative of matrix $M$ is determined by the evolution of its decomposed parts. For example, in the case of the singular value decomposition~\cite{uschmajew_geometric_2020},
\begin{align}
    \partial_t M &= \partial_t(U) \Sigma V^\dagger + U (\partial_t \Sigma)V^\dagger + U\Sigma \, \partial_t (V^\dagger) \nonumber
    \\
    &= \partial_t (U \Sigma) V^\dagger - U (\partial_t \Sigma)V^\dagger + U \partial_t (\Sigma V^\dagger)
    \label{eq:dt_M}
\end{align}
For the equation of motion $\partial_t M = F(M)$, the split-step procedure from time $t_0$ to time $t_1$, one solves 
\begin{align}
    & \partial_t (U \Sigma) = F(M) V(t_0) 
    &\rightarrow \tilde{M}_1 = U(t_1) \Sigma(t_1) V^\dagger(t_0) 
    \label{eq:te_us}
    \\
    & \partial_t \Sigma = U^\dagger(t_1) F(\tilde{M}_1) V(t_0) \, &\rightarrow \tilde{M}_0 = U(t_1) \tilde{\Sigma}(t_0) V^\dagger(t_0)
    \label{eq:te_s} 
    \\
    & \partial_t (\tilde{\Sigma} V^\dagger) = U^\dagger(t_1) F(\tilde{M}_0)
    & \rightarrow M(t_1)= U(t_1)\Sigma(t_1) V^\dagger(t_1)
    \label{eq:te_svt}
\end{align}
Note that the second equation is evolved backwards in time from $t_1$ to $t_0$. This algorithm is the DLRA with the projector splitting integrator (PS). 

In the case of dissipative systems, this backwards time evolution step can be unstable. This motivated the development of robust DLRA integrators, including the basis-update (BUG) method \cite{ceruti_robust_2024} and the alternating-projection (AP) method \cite{ ceruti_galerkin_2025}. In these methods, instead of performing backwards time evolution, the original state is projected onto the new basis $U(t_1)$, such that Eq.~\eqref{eq:te_s} is replaced by
\begin{align}
    & \tilde{\Sigma}(t_0) = U^\dagger(t_1) M(t_0)
    &\rightarrow \tilde{M}_0 = U(t_1) \tilde{\Sigma}(t_0) V^\dagger(t_0)
    \label{eq:ap_s} 
\end{align}
%

Extension to the tensor train case is straightforward, using the canonical decomposition $T=\textbf{L}M\textbf{R}$. However, to avoid performing the expensive computation of $\partial_t (\textbf{L}M)$ or $\partial_t (M\textbf{R})$, that step is separated into smaller steps by updating each of the tensor cores one at a time. Ultimately, one arrives at an iterative algorithm in which one performs the above steps at each tensor core of the tensor train.
At a high level, the DLRA-type procedure can be summarized as follows (see Sec.~\ref{sec:final} for more details). Starting with the center of orthogonality at the left-most tensor core, one performs the following at each tensor core:
\begin{enumerate}
    \item \textit{Blocking}. Construct the low-rank TT manifold $\left\{ |\textbf{L}^{(i)}_{\a_{i-1}}, \phi^{(i)}_{\s_i}, \textbf{R}^{(i)}_{\a_{i}} \rangle \right\}$, as well as the projector that projects a TT onto this manifold. The original equation of motion is projected onto that manifold to obtain an effective equation of motion. Provided that the TT is in canonical form, the manifold is defined by all tensor cores except for the active core (the center of orthogonality), and the effective equation of motion describes the equation of motion for the active core.
     
    \item \textit{Solving}. Solve the effective equation of motion (analogous to Eq.~\eqref{eq:te_us} or Eq.~\eqref{eq:te_svt}), obtaining the tensor core at the next time step. 
    
    \item \textit{Decimation}. If not the last tensor core to be updated, one uses the above solution to construct an updated low-rank manifold for the remaining tensor cores. This procedure will shift the center of orthogonality to the next site along the tensor train, which will be the active tensor core in the next iteration. 
    The standard projector-splitting DLR algorithm uses backwards time evolution (Eq.~\eqref{eq:te_s}), while the basis-update DLRA algorithm projects the original function onto the updated basis functions (Eq.~\eqref{eq:ap_s}.
    If it is the last tensor core to be updated, the tensor core is simply updated with the above result. At this point, the TT now represents the time evolved state. 
\end{enumerate}
(The terminology is taken from the density matrix renormalization group literature, e.g. in Ref.~\cite{Ronca_TDDMRG}).
In the DLRA literature, there is some variation in the ordering of the split step operations. However, this work will only consider sequentially sweeping through the tensor cores from left to right (steps are ordered as above) or right to left (the reverse). For the robust DLRA integrator, this is analogous to the AP method. 
%


\section{Interpolative DLRA}
\label{sec:interp}

This paper considers two ways of implementing DLRA for interpolative tensor trains. First is to modify the algorithm to be compatible with the interpolative TT representation. Second is to convert between the interpolative TT and the standard canonical format using oblique projections. In the remainder of the paper, the `X' will denote that the purely interpolative method is used, while `P' will denote that the method with oblique projections is used. The letter `G' will be used to denote the standard DLRA procedure that utilizes an orthogonal projection onto an orthonormal basis.

\subsection{Interpolative construction}
\label{sec:cur_dlra}

The structure of the DLRA algorithm remains the same. However, the underlying tensor core operations are modified by replacing the SVD operation with the CUR decomposition.




\subsubsection{Canonicalization and compression}
In the event that the original tensor is already in tensor train form, one uses the CUR decomposition to convert the tensor train into the cross-canonical form.

Assume that the center of orthogonality is at site $i$, meaning that the $i^\text{th}$ tensor core is a subselection of the original data, $M^{(i)}(\s_i)_{\a_{i-1},\a_{i}} = T(\mathcal{I}^{(i)}_{\a_{i-1}}, \s_i,\mathcal{J}^{(i)}_{\a_i})$. 
Once the CUR decomposition of $M^{(i)} = C^{(i)} \left(\hat{M}^{(i)}\right)^{-1} R^{(i)}$ is performed, the tensors are updated accordingly for left canonicalization (moving the center of orthogonality to site $i+1$):
\begin{align}
& M^{(i)} \leftarrow C^{(i)} \left(\hat{M}^{(i)}\right)^{-1} 
& M^{(i + 1)} \leftarrow R^{(i)} M^{(i+1)}
\end{align}
and for right canonicalization (moving the center of orthogonality to site $i-1$):
\begin{align}
& M^{(i)} \leftarrow \left(\hat{M}^{(i)}\right)^{-1} R^{(i)} 
& M^{(i - 1)} \leftarrow  M^{(i-1)} C^{(i)}
\end{align}
Note that instead of performing a CUR decomposition, one can instead do only row selection to obtain an XR decomposition for the left canonicalization procedure, or only do column selection to obtain a CX decomposition for the right canonicalization procedure.

As in the orthogonal representation, the canonicalization procedure does not alter the rank of the TT.  
%
To reduce TT rank, one needs to reduce the number of sub-selected rows (or columns). 
While one can use any scheme to perform the CUR decomposition with a predetermined rank, this work adaptively chooses the rank such that the approximation error is within some tolerance $\varepsilon$. Best results were obtained by first performing rank truncation using the SVD \cite{dolgov_parallel_2020}, and then the q-DEIM procedure is used on the reduced set of left singular vectors to obtain row selection indices $I$. The $(CU)R$ decomposition is performed on the original matrix. 
Note that because the SVD is used, and due to taking the inverse of the submatrix, there is no advantage in computational complexity of the interpolative construction compared to the orthogonal construction.

\subsubsection{Projection onto the low-rank manifold}


In the interpolative construction, the projection operator $\pmb{\Pi}$ evaluates the original data at the subselected indices, which are then used as the weights for the interpolating functions defined by $\textbf{L}$ and \textbf{R} (dropping the $(i)$ superscript for clarity):
\begin{align}
    \pmb{\Pi}  = \sum_{\a_{i-1},\a_{i}, \s_{i}} |\textbf{L}^{(i)}_{\a_{i-1}},\phi^{(i)}_{\s_i}, \textbf{R}^{(i)}_{\a_i} \rangle \langle {\mathcal{I}^{(i)}_{\a_{i-1}},\s_{i}, \mathcal{J}^{(i)}_{\a_i}}|
    \equiv \textbf{E} \, \textbf{S}
    \label{eq:cross_proj}
\end{align}
where $\textbf{E}$ is the interpolating functions from both $\textbf{L}^{(i)}$ and $\textbf{R}^{(i)}$, and $\textbf{S}$ is a matrix that selects the elements of $\s_{1:i-1}$ as specified by $\mathcal{I}^{(i)}$, the elements of $\s_{i+1:K}$ specified by $\mathcal{J}^{(i)}$, and the elements indexed by $\s_i$.
Applying the projector onto TT-vector,
\begin{align}
    \pmb{\Pi} \, T(\s_{1:K}) = T_{\mathcal{M}}(\s_{1:K}) = \sum_{\a_{i-1},\a_i} \textbf{L}^{(i)}_{\a_{i-1}}(\s_{1:i-1}) \,\, T \left[\,\mathcal{I}^{(i)}_{\a_{i-1}}, \, \s_i, \, \mathcal{J}^{(i)}_{\a_i} \right]
    \, \textbf{R}^{(i)}_{\a_{i}}(\s_{i+1:K}).
\end{align}
The notation $T[\mathcal{I},:,\mathcal{J}]$ denotes the appropriately selected elements of tensor $T(\s_{1:i-1}, \s_i, \s_{i+1:K})$.

These projectors are used to reduce a given problem from the full computational space down to the space of the TT manifold, defined by the tensor cores that are assumed to remain fixed. The TT operations and updates are thus limited to actions on a single tensor core. For example, matrix-vector multiplication is can be reduced to the tensor contraction of $A_\text{eff}$ with the single tensor core,
\begin{align}
    \pmb{\Pi} \, A \, \pmb{\Pi} \, x  = \textbf{E} \,  \underbrace{\textbf{S} \, A \, \textbf{E}}_{A_\text{eff}}  \, 
    \underbrace{\textbf{S} \, x }_{x_\text{eff}}
\end{align}
In most cases, $x$ already lives on the manifold of interest, i.e., $x = \textbf{L}M\textbf{R}$. If $x$ is in canonical form, then $x_\text{eff} = \textbf{S} \, x = x[\mathcal{I},:, \mathcal{J}]=M$. 
However, if $x$ is not in canonical form, which occurs if the rows and columns over-sampled, then one needs to explicitly compute $\textbf{S} x$.  

In the DLRA procedure, one must solve a reduced equation of motion at each tensor core. For example, consider the equation $\partial_t x = A f(x)$, where $A$ is a linear operator and $f$ is an elementwise operation. The reduced equation of motion on the TT manifold is then
\begin{align}
    \partial_t \pmb \Pi \, x = \textbf{E} \, \partial_t \underbrace{\textbf{S} \, x}_{x_\text{eff}} = 
    \underbrace{\textbf{S} A \, \textbf{E}}_{A_\text{eff}}  \, 
    \underbrace{ f(\textbf{S} \, x )}_{f(x_\text{eff})}
\end{align}

One can easily verify that the projection operation is valid, with $\pmb{\Pi}^2=\pmb{\Pi}$. 
%
%
If the tensor cores are in canonical form, satisfying Eqs.~\eqref{eq:canon_cx}, then $\textbf{S} \textbf{E}$ is the identity and one can readily show that
\[ \pmb{\Pi}^2 = \textbf{E} \,\underbrace{\textbf{S} \, \textbf{E}}_{\text{id.}} \, \textbf{S} = \textbf{E} \textbf{S} = \pmb{\Pi}. \]
If the tensor cores are not strictly canonical but the low-rank manifold is still low-rank in a least-squares sense, with $\textbf{E} = T (\textbf{S} T)^{-1} = T \, \hat{T}^{-1}$,
then
\[ 
\pmb{\Pi}^2 = (T \hat{T}^{-1}) \textbf{S} (T \hat{T}^{-1}) \textbf{S} =
T \underbrace{(\hat{T}^{-1}) \hat{T}}_{\text{id.}} \left(\hat{T}^{-1}\right) \textbf{S} = \pmb{\Pi}
\]
Because only the rows are oversampled in constructing submatrix $\hat{T}$, the submatrix is of full column rank and multiplying with the inverse on the left yields the identity.
A similar proof can be obtained for the projector projecting onto the right canonical tensor cores, and it immediately follows that Eq.~\eqref{eq:cross_proj} also is a valid projector.

\subsubsection{Subspace Expansion}
\label{sec:subspace_expansion_p}

A basis expansion scheme allows one to target multiple states with a single expanded manifold \cite{dolgov_polynomial_2015, dektor_interpolatory_2024}.
Consider a single tensor train with the center of orthogonality at tensor core $M^{(i)}_{\a_{i-1},\a_i}(\s_i)$.  The sizes of the $\alpha_{i-1},\alpha_{i}$ and $\s_i$ dimensions are $r_{i-1}, r_i$, and $d_i$, respectively. Also suppose there are multiple values of $M^{(i)}$ that the tensor core could be. Furthermore, one would like the TT manifold to contain all these states. 
Analogous to the method for orthonormal TTs, this is done by creating the expanded matrix $\tilde{M}$ by concatenating all matricized tensor cores $M^{(i),[k]}_{(\a_i,\s_i),\a_i}$ horizontally,
\begin{align*}
    \tilde{M}^{(i)} = 
    \begin{bmatrix} M^{(i),[1]}_{(\a_{i-1},\s_i),\a_i} & M^{(i),[2]}_{(\a_{i-1},\s_i),\a_i} & \hdots & M^{(i),[m]}_{(\a_{i-1},\s_i),\a_i}
    \end{bmatrix}
\end{align*}
and performing the CUR decomposition on this object to obtain an expanded set of interpolating functions $CU$ and row selection indices $I$. 
Before the expansion, the rank of the decomposition is limited to $\text{min}(r_{i-1} d_i, r_i)$. After expansion with $m$ tensors, the rank is bounded by $\text{min}(r_{i-1} d_i , m r_i)$, thus allowing one to sample with a larger number of points without oversampling.

Now with $m r_i$ indices selected, the TT is updated as follows
\begin{align}
    & M^{[k], (i+1)} \leftarrow  \sum_{\a_{i+1}} \left( M^{[k],(i)}_{(\a_{i-1},\s_i),\a_i)}[I,:] \right) M^{(i+1)}, 
    & M^{(i)} \leftarrow \tilde{M}[:,J] \left(\tilde{M}(I,J) \right)^{-1}.
\end{align}
The center of orthogonality is shifted right from $i$ to $i+1$, and there are now multiple values of the $(i+1)^\text{th}$ tensor core ($M^{(i+1),[k]}$). In contrast, there is only one value for the $i^\text{th}$ tensor core ($M^{(i)}$), which constructs the shared TT manifold.
Because the number of columns of $\tilde{M}$ is larger than each $\tilde{M}^{[k]}$, the subspace indexed by $\alpha_i$ has been expanded. 

In the context of DLRA AP-type algorithms, the different tensor cores could correspond to different stages from a multi-stage time integrator. Then, as the sweep across the tensor core proceeds, the TT will capture all stages of the time integrator. Of the $m$ tensor cores after the decimation step, one only needs to retain the tensor core corresponding to the initial value.

\subsection{Orthonormal construction with oblique projection}
\label{sec:oblique}

In this method, the tensor train is constructed in the orthonormal canonical form. However, one additionally constructs an oblique projector that converts between the orthonormal and interpolative representations so that one can readily obtain element-wise access to the original data. This method for tensor trains was introduced in Refs.~\cite{dektor_interpolatory_2024, dektor_collocation_2024}.
For clarity, this section uses the subscript $X$ to denote the interpolative basis and $Q$ to denote the orthonormal basis. For example,  $T = L_X M_X R_X$
denotes a TT in the interpolative basis, and 
$T = L_Q M_Q R_Q$
denotes the same TT in the orthonormal basis. 

Let us first consider the matrix case, with matrix $A=UM_Q$, where $U$ is a unitary and thus satisfies the left canonical condition. Obtaining the sub-sampled rows $M_X = A[I,:]$ is straightforward. Since $A[I,:] = U[I,:]M_Q$, one only needs to sample from $U$. 
The interpolating functionals are given by $U U[I,:]^{-1}$, as can be seen from the CUR decomposition of $A$. 
To convert back to the orthonormal representation, one simply inverts this expression, $M_Q = U[I,:]^{-1}M_X$. 

Extension to the TT case is straightforward. One can convert between the two representations with 
\begin{align}
    & M_X =  \textbf{L}_Q[\mathcal{I},:] \, M_Q \, \textbf{R}_Q[:, \mathcal{J}]
    \label{eq:Mq_to_Mx}
    \\
    & M_Q = \textbf{L}_Q[\mathcal{I},:]^{-1} M_X \textbf{R}_Q[:, \mathcal{J}]^{-1}
    \label{eq:Mx_to_Mq}
\end{align}
and the interpolative basis functions are
\begin{align}
    & \ket{\textbf{L}_X} = \ket{\textbf{L}_Q  \left(\textbf{L}_Q[\mathcal{I,:}]\right)^{-1} }
    \label{eq:q_to_x_left}
    \\
    & \ket{\textbf{R}_X} = \ket{ \left(\textbf{R}_Q [:,\mathcal{J}]\right)^{-1} \textbf{R}_Q}.
    \label{eq:q_to_x_right}
\end{align}
In the shortened notation introduced in the previous section, with $\textbf{E}_X$ denoting the interpolative basis and $\textbf{E}_Q$ denoting the orthonormal basis, and $\textbf{S}$ denoting the appropriate selection matrix,
\begin{align}
    & \textbf{E}_X = \textbf{E}_Q (\textbf{S} \, \textbf{E}_Q)^{-1} 
\end{align}
%
%
%
In contrast to the purely interpolative construction or purely orthonormal construction, one needs save both the orthonormal basis and the interpolating functions. Thus there is some additional computational overhead. 

\subsubsection{Canonicalization and compression}
The canonicalization and compression procedure is the same as the standard procedure, since the TT is still represented on an orthonormal basis. 
Assume that the center of orthogonality is at site $i$. 
First, the SVD decomposition is performed, $M_Q^{(i)} = U^{(i)} \Sigma^{(i)} {V^{(i)}}^\dagger$. The tensors are updated accordingly for left canonicalization (moving the center of orthogonality to site $i+1$):
\begin{align}
& M^{(i)}_Q \leftarrow U^{(i)} 
& M^{(i + 1)}_Q \leftarrow \Sigma^{(i)} R^{(i)} M^{(i+1)}
\end{align}
and for right canonicalization (moving the center of orthogonality to site $i-1$):
\begin{align}
& M^{(i)}_Q \leftarrow U^{(i)} \Sigma^{(i)} 
& M^{(i - 1)}_Q \leftarrow  {V^{(i)}}^\dagger
\end{align}
Note that instead of performing the SVD, one can instead perform the QR or LQ decomposition. 
However, the SVD is typically used if one intends to perform rank truncation.

However, once $U^{(i)}$ (or $V^{(i)}$) is obtained, an index selection procedure is performed on $U^{(i)}$ to obtain $I^{(i)}$ (or $J^{(i)}$). The matrices $U[I^{(i)},:]$ and its inverse are stored. 

\subsubsection{Projection onto the low-rank manifold}

The orthogonal projection operator onto the orthonormal basis is
\begin{align}
    \pmb{\Pi}_{QQ} = \sum |\textbf{L}^{(i)}_Q,\phi^{(i)},\textbf{R}^{(i)}_Q \rangle \langle \textbf{L}^{(i)}_Q,\phi^{(i)},\textbf{R}^{(i)}_Q |  = \textbf{E}_Q \textbf{E}_{Q}^\dagger
\end{align}
and the projection operator onto the interpolative basis is given in Eq.~\eqref{eq:cross_proj}.
%
The oblique projectors that projects a TT from the interpolative basis to the orthonormal basis, and vice verse, are 
\begin{align}
    & \pmb{\Pi}_{QX} = \textbf{E}_Q \textbf{E}_Q^\dagger  \textbf{E}_X \textbf{S} = \textbf{E}_Q (\textbf{S}\textbf{E}_Q)^{-1} \textbf{S}
    \\
    & \pmb{\Pi}_{XQ} = \textbf{E}_X \textbf{S} \textbf{E}_Q \textbf{E}_Q^\dagger,
\end{align}
which give rise to Eqs.~\eqref{eq:Mx_to_Mq} and ~\eqref{eq:Mq_to_Mx} when the two projectors are applied to a tensor train of the interpolative and orthonormal canonical form, respectively. 

Again, consider the dynamics of $\partial_t x = A f(x)$ on the tensor train manifold, where $A$ is a linear operator and $f$ is some arbitrary element-wise operator. Depending on the TT manifold of interest, one can arrive at different equations of motion.
%
%
For the orthonormal TT manifold, one can utilize the interpolative representation to perform the nonlinear operations and then project the equation of motion back onto orthonormal basis:
\begin{align}
    \partial_t \pmb{\Pi}_{QQ} \, x = \textbf{E}_Q \,  \partial_t \textbf{E}_Q^\dagger  \, x 
    & = \textbf{E}_Q \underbrace{\textbf{E}^\dagger_Q A \textbf{E}_Q }_{A_\text{eff}} 
    \, \underbrace{\textbf{E}_Q^\dagger \textbf{E}_X f(\textbf{S} \textbf{E}_Q \textbf{E}_Q^\dagger x)}_{\pmb{\Pi}_{QX} f(S x)}
    \nonumber 
    \\[0.2cm]
    & = \textbf{E}_Q \, \underbrace{\textbf{E}^\dagger_Q A \textbf{E}_Q }_{A_\text{eff}} 
    \, \underbrace{(\textbf{S} \textbf{E}_Q)^{-1} f(\textbf{S} \textbf{E}_Q \textbf{E}_Q^\dagger x)}_{\pmb{\Pi}_{QX} f(Sx)} \, .
    \nonumber
\end{align}
However, the quality of the approximation depends on how well $E_Q$ can capture the nonlinear contribution $f(x)$ in addition to $x$.
In the event that $f$ is linear, the $f(x_\text{eff})$ term can be greatly simplified. 

The equivalent dynamics on the interpolative manifold is given by
\begin{align}
    \partial_t \pmb{\Pi}_{XX} \, x = \textbf{E}_X \,  \partial_t \textbf{S} \textbf{E}_Q \textbf{E}_Q^\dagger  \, x 
    = & \textbf{E}_X \underbrace{\textbf{S} \textbf{E}_Q \textbf{E}^\dagger_Q A \textbf{E}_Q (\textbf{S} \textbf{E}_Q)^{-1} }_{A_\text{eff}} 
    \, \underbrace{f(\textbf{S} \textbf{E}_Q \textbf{E}_Q^\dagger x)}_{f(S x)}
    \nonumber
    \\
    = & \textbf{E}_X \, \underbrace{\textbf{S} \textbf{E}_Q \textbf{E}^\dagger_Q A \textbf{E}_X }_{A_\text{eff}} 
    \, \underbrace{f(\textbf{S} \textbf{E}_Q \textbf{E}_Q^\dagger x)}_{f(Sx)} \, .
    \nonumber
\end{align}
One can verify that they are equivalent by noting that 
\[ \partial_t \textbf{E}^\dagger_Q x = (\textbf{S} \textbf{E}_Q)^{-1} \partial_t \textbf{S} \textbf{E}_Q \textbf{E}_Q^\dagger x  = (\textbf{S} \textbf{E}_Q)^{-1} \textbf{S} \textbf{E}_Q \textbf{E}_Q^\dagger A = \textbf{E}_Q^\dagger A \] if $(\textbf{S} \textbf{E}_Q)$ is invertible. This is generally the case if no oversampling is performed ($r$ indices are selected for a rank-$r$ matrix), and if the selection indices are properly selected, for example, via the maxvol algorithm or DEIM.

Note that the equivalence is due to the projection onto the low-rank orthonormal basis, 
$\textbf{E}_Q \textbf{E}_Q^\dagger$. Without this projection, the equations of motion are no longer equivalent, since 
\( (\textbf{S} \textbf{E}_Q)^{-1} \textbf{S} A \neq \textbf{E}_Q^\dagger A \) because $\textbf{S}$ is not left invertible.

\subsubsection{Subspace Expansion}
\label{sec:subspace_expansion}

The basis expansion procedure of targeting multiple tensor cores $M^{(i),[k]}$ follows the standard procedure, and is analogous to that described earlier for the interpolative TT construction. 
If performing a left canonicalization procedure, one constructs the expanded matrix $\tilde{M}$ by concatenating all matricized tensor cores $(M^{(i),[k]}_Q)_{(\a_i,\s_i),\a_i}$ horizontally.
One then performs the QR or SVD decomposition, and the left unitary is the new expanded orthonormal basis. A row-index selection procedure is used on $U$ to determine $I^{(i)}$.

However, one must first project all of the tensor cores onto the orthonormal basis using Eq.~\eqref{eq:Mx_to_Mq} if they are not already in that space.

\subsection{Final interpolative dynamical low-rank algorithm}
\label{sec:final}

\begin{figure}
    \centering
    \includegraphics[width=0.95\linewidth]{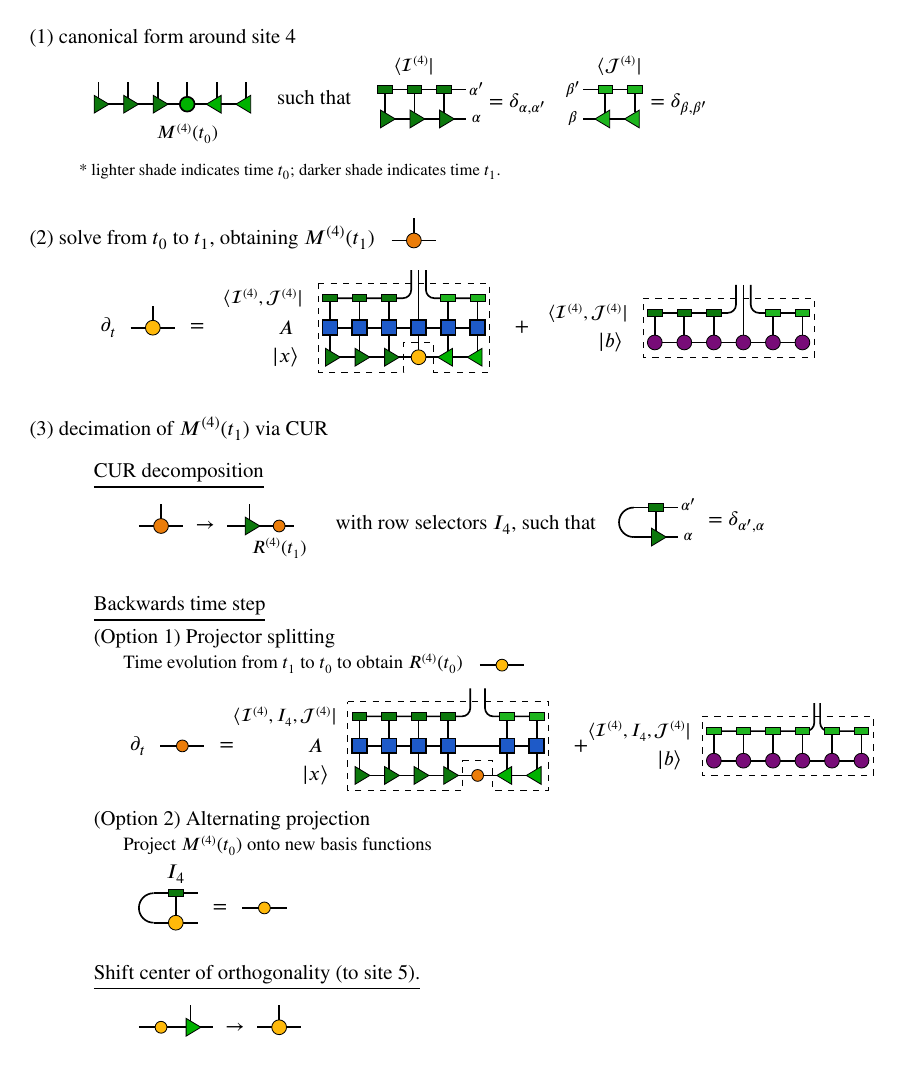}
    \caption{An overview of the three steps in a single iteration of the dynamical low-rank time integration algorithm (see Sec.~\ref{sec:final}) with the interpolative TT construction (DLR-X), solving $\partial_t x + Az = b$. This figure explicitly shows the DLR procedure with the center of orthogonality at the fourth tensor core, with the sweeping procedure going from left to right. However, the procedure is iterated over all cores in the tensor train. In the diagrams, the TT with green cores depicts the vector $x$, the TT with blue cores depicts the operator $A$, and the TT with the purple cores depicts the vector $b$. The right/left-pointing triangles denote left/right canonical tensors in the interpolative construction, and the flat rectangles denote index selection.}
    \label{fig:dlrx}
\end{figure}

\begin{figure}
    \centering
    \includegraphics[width=0.93\linewidth]{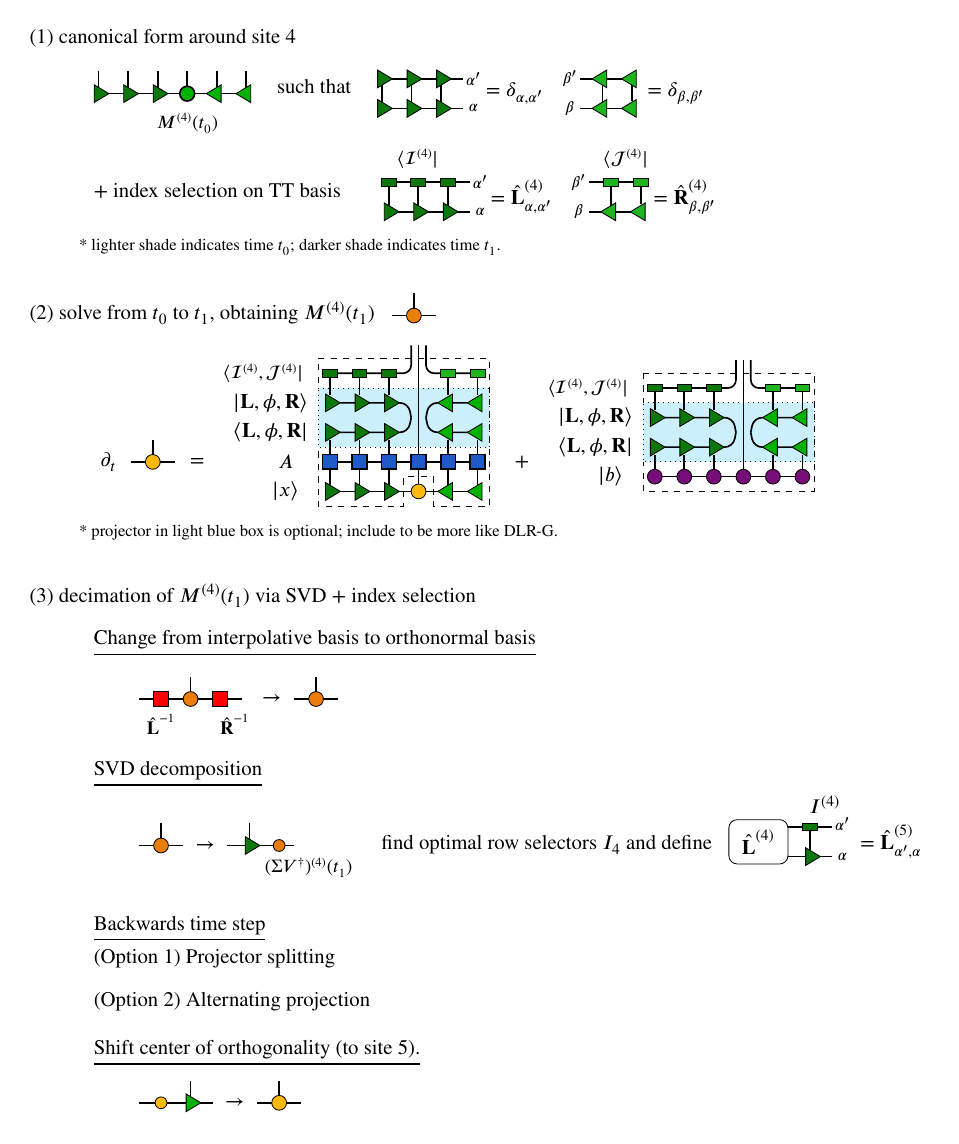}
    \caption{An overview of the three steps in a single iteration of the dynamical low-rank time integration algorithm (see Sec.~\ref{sec:final}) with oblique projection (DLR-P), solving $\partial_t x + Ax = b$. This figure explicitly shows the DLR procedure with the center of orthogonality at the fourth tensor core, with the sweeping procedure going from left to right. However, the procedure is iterated over all cores in the tensor train. Here, the right/left-pointing triangles denote left/right canonical tensors in the orthonormal construction.}
    \label{fig:dlrp}
\end{figure}

Consider the time-dependent PDE of the form $\partial_t u + A F(u) = b$.
The DLR-type procedure for the interpolative construction and with the oblique projector (depicted diagrammatically in Fig.~\ref{fig:dlrx} and Fig.~\ref{fig:dlrp}, respectively) is as follows:

\begin{enumerate}
    \item Canonicalize the TT such that the center of orthogonality is at the leftmost tensor core.
    \item Expand the TT basis as desired to improve the accuracy of the calculation.
    \begin{itemize}
        \item If there is a non-zero source term $b$, then the manifold of $u$ should be expanded such that both $b$ and $u$ can be represented on the TT manifold up to the desired accuracy.
    \end{itemize}
    \item Build the basis of the time evolved state iteratively, starting from the leftmost tensor core and then sweeping from left to right. At each iteration, only the $i^\text{th}$ tensor core is modified, while the others remain fixed. The TT must be in canonical form, with the center of orthogonality at site $i$. At each iteration, the following procedure is performed:
    
    \begin{enumerate}
        \item \textit{Blocking}. Construct the low-rank manifold defined by all tensor cores except for the $i^\text{th}$ core, 
        \[ \mathcal{M}^{(i)}  = 
        \left\{ \ket{ \lenv{i}, 
        \phi^{(i)}_{\s_i}, 
        \renv{i} 
        } \right\}
        \quad \forall \, \a_{i-1} \in [1, r_{i-1}], \, \a_{i} \in [1, r_i], \, \s_i \in \{0,1\} \]
        as well as the projector $\pmb{\Pi}$ that projects the TT onto the interpolative manifold. Using this projector, obtain an effective equation of motion for active tensor core. For example, if the equation of motion is $\ddt u + A F(u) = b$, then the effective equation of motion is
        \begin{align}
            \partial_t M^{(i)}_X + A_\text{eff} F(M_X^{(i)}) = b_\text{eff},
        \end{align}
        where, for interpolative DLRA,
        \begin{align}
            b^\text{eff}_{(\a_{i-1}, \s_i, \a_i)}
            &= b \left[ \isel{i}, \s_i, \jsel{i} \right] \\
            A^\text{eff}_{(\a'_{i-1},\s'_i,\a'_i),(\a_{i-1},\s_i,\a_i)} 
            &= \left( A \,  \ket{\lenv{i}, \phi_{\s_i}, \renv{i}} \right) \left[ \bisel{i}, \s'_i, \bjsel{i} \right]
        \end{align}
        and for DLRA with oblique projection,
        \begin{align}
            b^\text{eff}_{(\a_{i-1}, \s_i, \a_i)}
            &= (\Pi_Q b) \left[ \isel{i}, \s_i, \jsel{i} \right] \\
            A^\text{eff}_{(\a'_{i-1},\s'_i,\a'_i),(\a_{i-1},\s_i,\a_i)} 
            &= \left( (\Pi_Q A)  \,  \ket{\lenv{i}, \phi_{\s_i}, \renv{i}}_X \right) \left[ \bisel{i}, \s'_i, \bjsel{i} \right]
        \end{align}
        where Q and X are again used to denote if the projection and basis vectors correspond to the orthonormal or the interpolative space. 
         
        \item \textit{Solving}. Solve the effective equation of motion, obtaining the tensor core at the next time step, $M^{(i)}_X(t+\dt)$. In this work, explicit Euler, fourth-order Runge-Kutta (RK4), and the Crank-Nicolson time integrators are considered.

        For DLRA with oblique projection, the resulting $M^{(i)}_X(t+\dt)$ need to be projected back onto the orthonormal TT manifold, obtaining $M^{(i)}_Q(t+\dt)$.
        Assuming the tensor cores are in the desired space, the subscripts 'Q' and 'X' will now be dropped.
        
        \item \textit{Decimation}. 

        If at the last tensor core, update it with $M^{(i)}(t+\dt) $.
        
        If not at the last tensor core,
        decompose the updated tensor: $M^{(i)}(t+\dt) = \bar{U}^{(i)} W$, where $\bar{U}^{(i)}$ is in left canonical form and is used to construct an updated manifold for the $(i+1)^\text{th}$ tensor core. One also needs to determine the row selection indices $I^{(i)}$. If rank adaptivity is allowed, the decimation procedure is performed with error threshold $\varepsilon_{in}$.
        \begin{itemize}
            \item In the projector splitting (PS) algorithm, one performs backwards time evolution on the weight matrix $W$ before contracting it with the next tensor to shift the center of orthogonality.
            \begin{enumerate}
                \item Update tensor core $i$: $M^{(i)} \leftarrow \bar{U}^{(i)}$
                \item Update the left environment to include $\bar{U}^{(i)}$, yielding $\textbf{L}_{\a'_i}^{(i+1)}$
                \item Backward time evolution of $W^{(i)}(t+\dt)$ along the manifold $\{ |\textbf{L}_{\a'_i}^{(i+1)}, \renv{i} \rangle \}$, obtaining $W^{(i)}(t)$
                \item Update $M^{(i+1)} \leftarrow W^{(i)} M^{(i+1)}$; the center of orthogonality is now at site $i+1$. 
                \item Obtain the right environment without $V^{(i+1)}$, which is $\textbf{R}_{\a'_{i+1}}^{(i+1)}$
            \end{enumerate}

            In the single-site PS algorithm presented here, the rank of the TT remains unchanged. To enable rank-adaptivity, one can use a two-site variant, in which one performs forward time evolution on a supercore consisting of two sites ($i$ and $i+1$) during the ``solve'' step, and then backwards time evolution on one site ($i+1$) during the ``decimation'' step.

            \item
            In the alternating-projection (AP) DLR algorithms, one projects the original tensor $M^{(i)}$ onto the new basis. This is done by updating the sites as
            \begin{enumerate}
                \item $M^{(i)} \leftarrow \bar{U}^{(i)}(t+\dt) $
                \item $M^{(i+1)} \leftarrow \sum_{\a'_i, \a'_{i-1}, \s_i} \left( \bar{S}^{(i)}_{\a'_{i-1},(\a_i,\s_i)} \right)^* M^{(i)}_{\a'_{i-1},\a'_i}(\s_i) V_{\a'_{i},\a_{i+1}}^{(i+1)}(\s_{i+1}) $; the center of orthogonality is now at site $i + 1$
            \end{enumerate}
            Targeting the initial state helps ensure that the TT remains unchanged as one sweeps through all the tensor cores. Instead, the sweeping procedure iteratively builds a new TT manifold to accurately capture the next time step.
            
            For higher order time integration, instead of using the $\bar{U}^{(i)}$ defined above, one expands the basis by targeting multiple states. This allows one to obtain an enlarged basis $\bar{A}(i)$ that targets all of the intermediate stages. 
            For Euler time-stepping, though not strictly necessary for first-order accuracy, one constructs the canonical tensor $\bar{U}^{(i)}$ by targeting the time evolved state as well as the initial state, \[ \begin{bmatrix} M^{(i)}(t) &  M^{(i)}(t+\dt) \end{bmatrix} = \bar{U}^{(i)} W. \] 

            For the CN scheme, the same tensor expansion is used
            \[ \begin{bmatrix} M^{(i)}(t) & M^{(i)}(t+\dt) \end{bmatrix} \]

            In RK4, one also targets the state at additional intermediate time steps, targeting the expanded matrix (see Eq. (5) in \cite{feiguin2005time}),
            \[ \begin{bmatrix} M^{(i)}(t) & M^{(i)}(t + \frac{\dt}{3}) & M^{(i)}(t + \frac{2\dt}{3}) & M^{(i)}(t+\dt) \end{bmatrix}. \]
            One can achieve similar results by additionally targeting the derivatives at each stage \cite{nobile_robust_2025}. For example, if $k_m$ is derivatives at the $m^\text{th}$ stage, one could also consider
            \[ \begin{bmatrix} M^{(i)}(t) & k_1 & k_2 & k_3 & M^{(i)}(t+\dt) \end{bmatrix}. \]

            For explicit time integrators, instead of basis expansion, another option is to save separate TT manifolds for each of the intermediate states \cite{carrel_interpolatory_2025}. This would allow one to use lower ranks for each TT, but one would need to compute the effective equation of motion for each. However, this could be done in parallel, and may be an interesting avenue for performant codes. However, this is not done in this work.
            
        \end{itemize}

    \end{enumerate}
    
    \item With the update of the last tensor core, the TT has finally been moved forward in time. The rank of the tensor train at this stage is denoted as $r_{in}$

    \item Rank truncation procedure (sweeping in the opposite direction) to the desired rank or truncation threshold $\varepsilon$ if necessary. The rank of the tensor train at this stage is denoted as $r$.

\end{enumerate}
Better results are obtained if the internal threshold $\varepsilon_{in}$ is smaller than the truncation threshold $\varepsilon$ used for the final rank truncation procedure. However, unless otherwise specified, $\varepsilon_{in}=\varepsilon$ is used.
Currently, the final rank truncation is performed using the SVD truncation procedure since some claim better stability \cite{zheng_SLAR_2024}.
However, if the TT is in the interpolative basis, CUR-based rank truncation can also be used, which would avoid the need to first convert the QTT into the orthonormal canonical form.

\clearpage

\subsubsection{Two-Site Variant}
During time integration, the rank of the true solution may increase. As such, the DLR algorithm should allow for the tensor train rank to adapt accordingly. 
In the AP method, rank adaptivity is achieved via basis expansion, as described above. However, this is not done in the PS method. Previous works achieved rank adaptivity by oversampling the matrix. In the orthogonal construction, this is done by augmenting the decomposition of the matrix with unit basis vectors or random orthonormal vectors that are multiplied by zero \cite{hochbruck_rank-adaptive_2023, lubich_projector-splitting_2014}. In the interpolative construction, this was achieved by using multiple index selection methods to increase the number of rows or columns selected, oversampling the matrix in the CUR decomposition \cite{dektor_collocation_2024}. 

However, a more systematic method for rank adaptivity is to use a two-site variant. In this case, the tensor train is of the form 
\begin{align}
    T(\s_{1:K}) = \sum_{\a_{i-1},\a_{i}} \textbf{L}_{\a_{i-1}} (\s_{1:i-1}) \, M^{(i)}_{\a_{i-1},\a_{i}} (\s_i) \, M^{(i+1)}_{\a_{i},\a_{i+1}} (\s_{i+1})\, \textbf{R}_{\a_{i+1}} (\s_{i+2:K}) \, .
\end{align}
such that the center of orthogonality is at two adjacent tensor cores. By contracting the two tensors together, generating a supercore $\textbf{M}(\sigma_i,\sigma_{i+1}) = \sum_{\alpha_{i+1}} M^{(i)}_{\alpha_{i-1},\alpha_i}(\sigma_i) M^{(i+1)}_{\alpha_{i},\alpha_{i+1}} (\sigma_{i+1}) $, the DLRA algorithm is the same as before.

Consider the equation of motion $\partial_t x = F(x)$.
If sweeping left to right, the DLRA procedure for a single iteration is
\begin{enumerate}
    \item \textit{Blocking}. Already satisfied.
    \item \textit{Solving}. Solve from time $t_0$ to $t_1$ to obtain $\textbf{M}(t_1)$.
    %
    \item \textit{Decimation}. Separate the two sites on the supercore (with error threshold $\varepsilon_{in}$),
    \[ \textbf{M}_{(\a_{i-1},\s_i);(\s_{i+1},\a_{i+1})}^\textbf{1}(t_1) = \sum_{\a_i} U^{(i)}_{(\a_{i-1},\s_i),\a_i}(t_1) \tilde{M}^{(i+1)}_{\a_i;(\s_{i+1},\a_{i+1})}(t_1), \] 
    \[ \tilde{T}_1 = \underbrace{\textbf{L}^{(i-1)} U^{(i)}( t_1)}_{\textbf{L}^{(i)}} \tilde{M}^{(i+1)}(t_1) \textbf{R}^{(i+2)} \]
    Obtain $\tilde{M}^{i+1}(t_0)$ by performing backwards time evolution from time $t_1$ to $t_0$.
    %
    Generate the new supercore for the next iteration, $\textbf{M}(\s_{i+1},\s_{i+2}) =M^{(i+1)}(\s_{i+1}) B^{(i+2)}(\s_{i+2})$.
\end{enumerate}
For simplicity, let all ranks in the QTT be equal to $r_{in}$, such that the tensor core $M^{(i)}_{(\a_{i-1},\s_i), \a_{i+1}}$ is of size $(r_{in} \times d) \times r_{in}$. 
In the single-site case, the maximum rank of this matricized tensor is $r_{in}$. As a result, during the decimation step, the rank of the tensor train is not increased.
In the two-site case, the size of the supercore $\textbf{M}_{(\a_{i-1},\s_i), (\a_{i+1},\s_{i+1})}$ is $(r_{in} \times d) \times (r_{in} \times d)$. The maximum rank of the matricized supercore is $dr_{in}$. Thus, the size of the basis at tensor core $i$ can be increased from $r_{in}$ to $dr_{in}$.

For the two-site algorithm, denoting the rank of the QTT-operator $A$ as $r_A$, the cost of tensor contraction is roughly $\mathcal{O}(r_{in}^3 r_A d^2 + r_{in}^2 r_A^2 d^2)$ and the cost of canonicalization is $\mathcal{O}(r_{in}^3 d^3)$. 
In comparison, for the one-site case variant, the cost of the tensor contraction is roughly $\mathcal{O}(r^3_{in} r_A d + r^2_{in} r_A^2 d)$ and the cost of canonicalization is $\mathcal{O}(r^3_{in} d^3)$.
Therefore, for the QTT construction where $d=2$, the two-site variant is not prohibitively more expensive than the one-site variant.


\section{Numerical Examples}

\subsection{Inviscid Burger's equation}

Let us consider the inviscid Burger's equation in 1-D in flux form,
\[ \partial_t u + \frac{1}{2} \partial_x u^2 = 0,\]
where $u$ is the velocity field that varies in space (along $x$) and over time. 
In the following simulations, the spatial domain $[0,1)$ is uniformly discretized with $2^L$ grid points, such that $u$ can be represented as a QTT of length $L$ with each quantized dimension of size $d=2$. A time step size of $\Delta t =0.9 \Delta x$ is used, where $\Delta x$ is the grid spacing.

The inviscid Burger's equation is often solved with an upwind time integration scheme. 
Let $u_j^n$ denote the discretized velocity field with subscript $j$ denoting the position in space and superscript $n$ denoting its position in time.
The velocity field $u$ at the next time step is given by   
\begin{align}
    u_j^{n+1} = u^n_j - \frac{\Delta t}{\Delta x} \left( \left(\tilde{u}_{j+1/2} \right)^2 - \left( \tilde{u}_{j-1/2} \right)^2 \right)
    \label{eq:upwind_euler_burgers}
\end{align}
where 
\begin{align}
    \tilde{u}_{j+1/2} =
    \begin{dcases}
        \quad u_{j+1} & \text{if } u_{j+1} > u_j \text{ and } \frac{1}{2} \left( u_{j+1} + u_j \right) > 0 \\
        \quad u_{j} & \text{if } u_{j+1} > u_j \text{ and } \frac{1}{2} \left( u_{j+1} + u_j \right) > 0 \\
        \quad u_{j+1} & \text{if } u_{j+1} < u_j \text{ and } u^+ > 0 \\
        \quad u_{j} & \text{if } u_{j+1} < u_j \text{ and } u^- < 0 \\
        \quad  0 & \text{o.w.}
    \end{dcases}
    \label{eq:upwind_burgers}
\end{align}

Now consider performing upwind time evolution using DLRA. 
With the orthonormal construction, when $u$ is in canonical form, it is not possible to compute the desired flux at each grid point because since the center of orthogonality does not directly correspond to the values of $u$ at any specific grid points.
In contrast, with the interpolative construction, the center of orthogonality is exactly the values of $u^n$ at the grid points determined by selection indices $\mathcal{I}$ and $\mathcal{J}$. The values of $u^n$ at the adjacent grid points can be measured efficiently, and then one can compute the flux at selected points. 

The subroutine in the DLR-X procedure at each tensor core is as follows:
\begin{enumerate}
    \item blocking: determine the manifold from frozen tensor cores of $u$. Let $S^+$ and $S^-$ denote the operator that shifts all entries left and right by one, respectively. Project the operators onto the low-rank manifold,
    \begin{align}
        (S^{\pm}_{\text{eff}})_{\a'_{i-1},\s'_i,\a'_i;\a_{i-1},\s_i,
        \a_i} = (S^{\pm} |\lenv{i}, \phi_{\s_i}, \renv{i}\rangle )[\mathcal{I}^{(i)}_{\a'_{i-1}}, \s'_i, \mathcal{J}^{(i)}_{\a'_i}]
    \end{align}
    \item solving: Let $\mathcal{I}$ denote the indexed grid points. Determine
    \[ u^n[\mathcal{I}] = M^{(i)} \]
    \[ u^n[\mathcal{I}\pm1] = S^{\pm}_\text{eff} M^{(i)} = (S^{\pm} u^n)_{\mathcal{I}} = u^n_{\mathcal{I}\pm 1} \]
    Then, solve for the fluxes at each index $\mathcal{I}$ and compute the field at the next time step using Eq.~\eqref{eq:upwind_euler_burgers}.
    \item decimation: shift the center of orthogonality to the next tensor.
    \begin{itemize}
        \item Alternating projection (AP): expand the subspace of $u$ to target both $u^n$ and $u^{n+1}$, following the procedure described in Sec.~\ref{sec:final}: perform CUR decomposition of the matrix $\tilde{M} = [u^{n+1}_{(\a_{i-1},\s_i),\a_i}, \, u^n_{(\a_{i-1},\s_i),\a_i}] = CUR$. Update site $i$ with the interpolating functions $A^{(i)}_{\a_{i-1}, \a_i}(\s_i) =(CU)$, and retain the row selection indices $I^{(i)}$. Project site $i+1$ of $u$ onto this new set of interpolating functions, following $M^{(i+1)} \leftarrow RM^{(i+1)}$.
        \item Projector splitting (PS): This scheme is not well suited for this problem. In theory, one can also use an upwind procedure (taking proper care of the sign change in the time step).
        However, doing so adds additional dissipation.
    \end{itemize}
\end{enumerate}
The subroutine for DLR-P follows the same overall steps.

\begin{figure}
    {\small
    \hspace{0.8cm} (a) shock formation \hspace{0.2cm} (b) shock propagation \hspace{0.3cm} (c) rarefaction }
    \\
    \includegraphics[width=\linewidth]{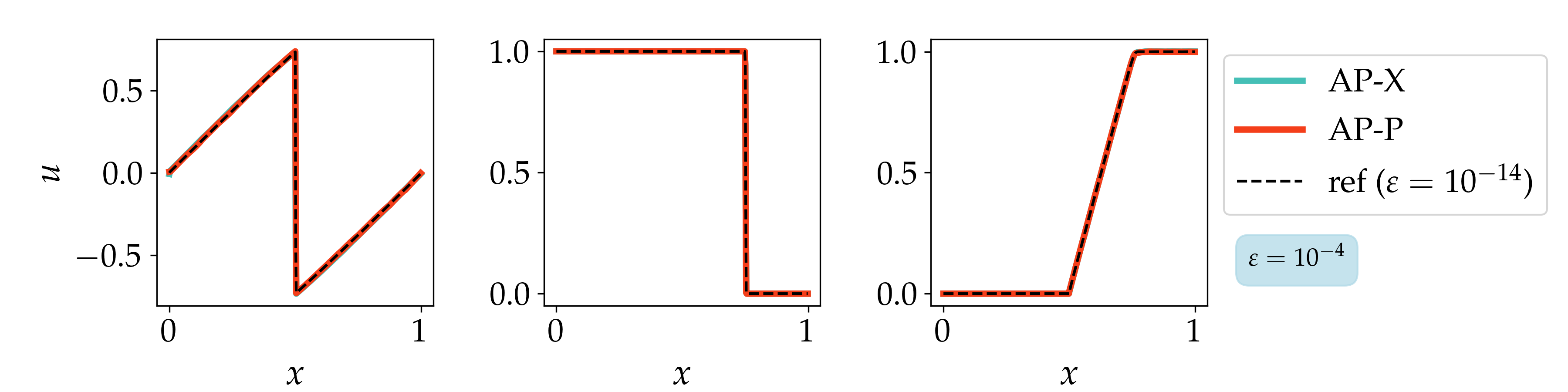}
    \includegraphics[width=\linewidth]{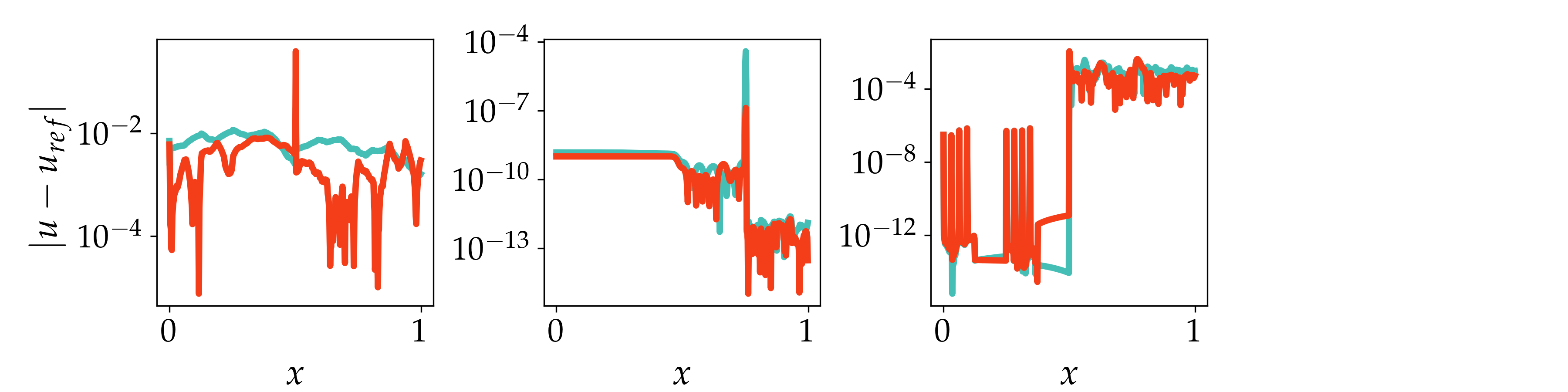}
    \includegraphics[width=\linewidth]{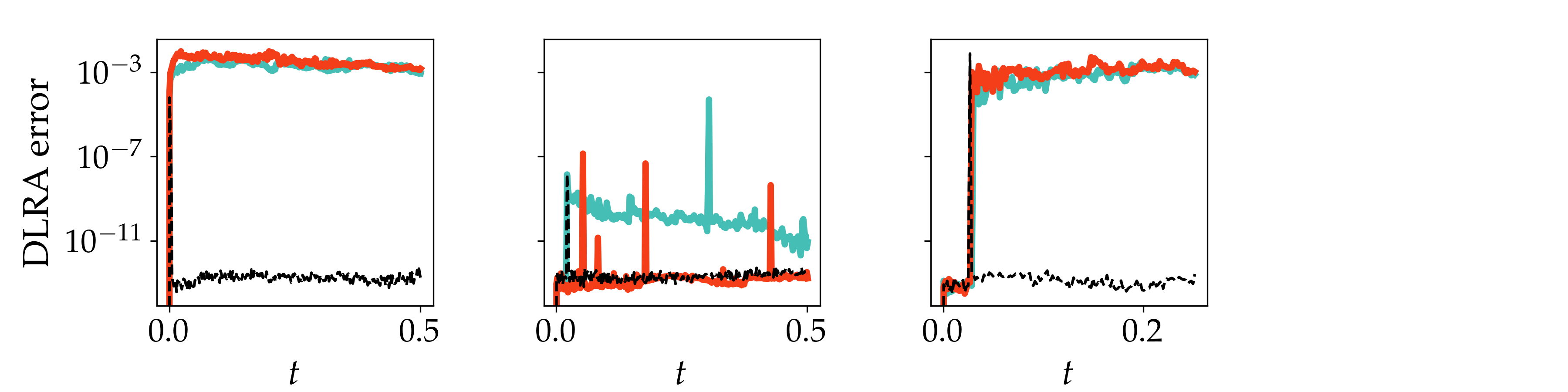}
    \includegraphics[width=\linewidth]{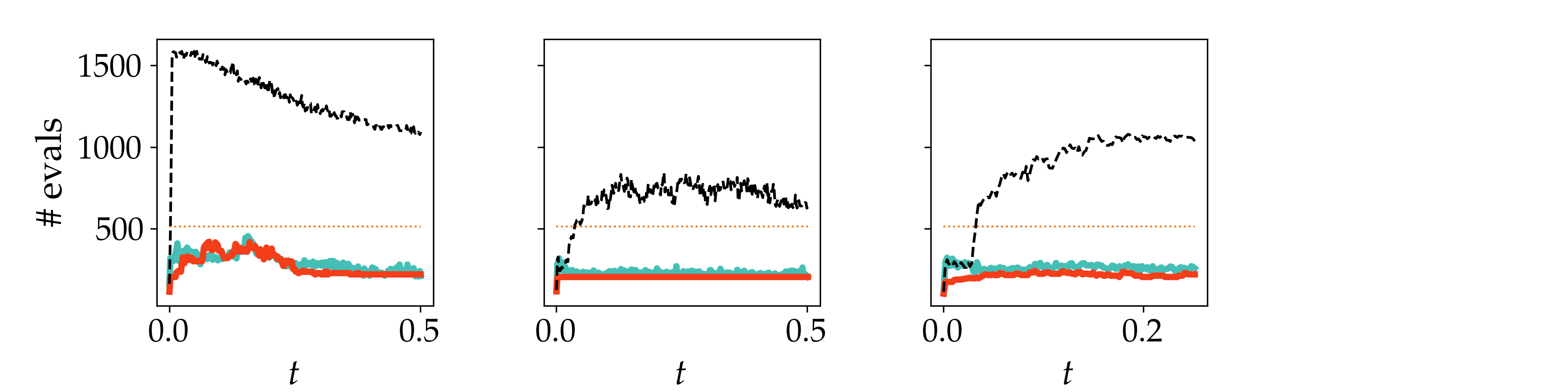}
    \caption{DLR time evolution with rank truncation threshold $\varepsilon=10^{-4}$ on a uniform grid with $2^9$ grid points. Test problems include (a) shock formation, with initial condition $\sin(2\pi x)$, (b) shock propagation, with initial condition $1-H(x-0.5)$, where $H(x)$ is the Heaviside step function, and (c) rarefaction, with initial condition $H(x-0.5)$. The DLR-X calculations are performed with each of the specified truncation thresholds. (Top) the time-evolved velocity field at $T=0.5$ for (a) and (b), and $T=0.25$ for (c). (Second row) the error of field $u$ with respect to the $\varepsilon=10^{-14}$ result, denoted as the reference. 
    (Third row) the error arising from the DLR approximation over the course of the simulation. 
    (Bottom) the total number of evaluations of Eq.~\eqref{eq:upwind_euler_burgers}. For the dense calculation, the total number of evaluations is $2^9$, and is depicted by the thin dotted line.}
    \label{fig:burgers}
\end{figure}

Results for calculations with $L=9$ and various initial conditions are shown in Fig.~\ref{fig:burgers}. A moderate error threshold of $\varepsilon=\varepsilon_{in}=10^{-4}$ is used. The calculation is performed with the purely interpolative scheme.
The first row depicts the velocity profile at the final time. Only at very aggressive tolerances are the errors visible by eye. The second row depicts the error of the calculation with respect to a $\varepsilon=10^{-14}$ calculation. The third row depicts the DLRA error at each time step, defined as the $L_2$ difference between the $u^{n+1}$ obtained from the DLRA procedure before the final compression step and that obtained from an exact upwinding procedure when starting from the same initial state $u^n$. This measures the error that occurs from projecting the dynamics onto the reduced TT manifold. The fourth row depicts the number of Euler evaluations performed per time step. This value is bounded by $L r_{in}^2 d$, where $L=9$ and $d=2$ for these calculations. 
Due to the relatively low-rank structure of the solution, the AP-DLRA time integration scheme allows for fewer evaluations compared to the dense calculation, which requires 512 evaluations at each time step. However, the number of evaluations does not directly reflect the computational cost, since one must also consider the additional overhead of QTT algorithms due to the tensor contractions and decompositions, whose computational costs scale like $\mathcal{O}(d^2 r_{in}^2 r)$, where usually $r < r_{in}$. However, for a given problem, the ranks $r$ and $r_{in}$ as well as the quantization size $d$ are expected to be independent with respect to grid resolution.
Again, the QTT method is only expected to possibly outperform traditional numerical methods for calculations requiring high resolution, as well as those of higher dimensionality.

For this example, DLR-X and DLR-P yield comparable results, with DLR-P allowing for slightly smaller ranks and slightly smaller errors. This is consistent with the fact that orthonormal TT construction is optimal in the $L_2$ sense, whereas the interpolative TT construction is not.

In some cases, the low-rank calculations were susceptible to issues due to the rank being too small, resulting in inaccurate DLRA time integration. The shock propagation test suffered from this issue the most. On occasion, the shock front would not be propagated forward, thus leading to incorrect shock propagation speeds. As such, it was helpful to enforce a minimum rank. In these calculations, the minimum rank is set to four.

\subsection{Maxwell's equation with non-uniform dielectric}

In this next example, we consider a 2-D simulation of an electromagnetic wave trapped inside a dielectric cavity, though the cavity is smaller than the simulation domain size.

\begin{figure}
    \centering
    \begin{subfigure}{\linewidth}
        \includegraphics[width=\linewidth, trim={0cm 0.2cm 0cm 0cm}, clip]{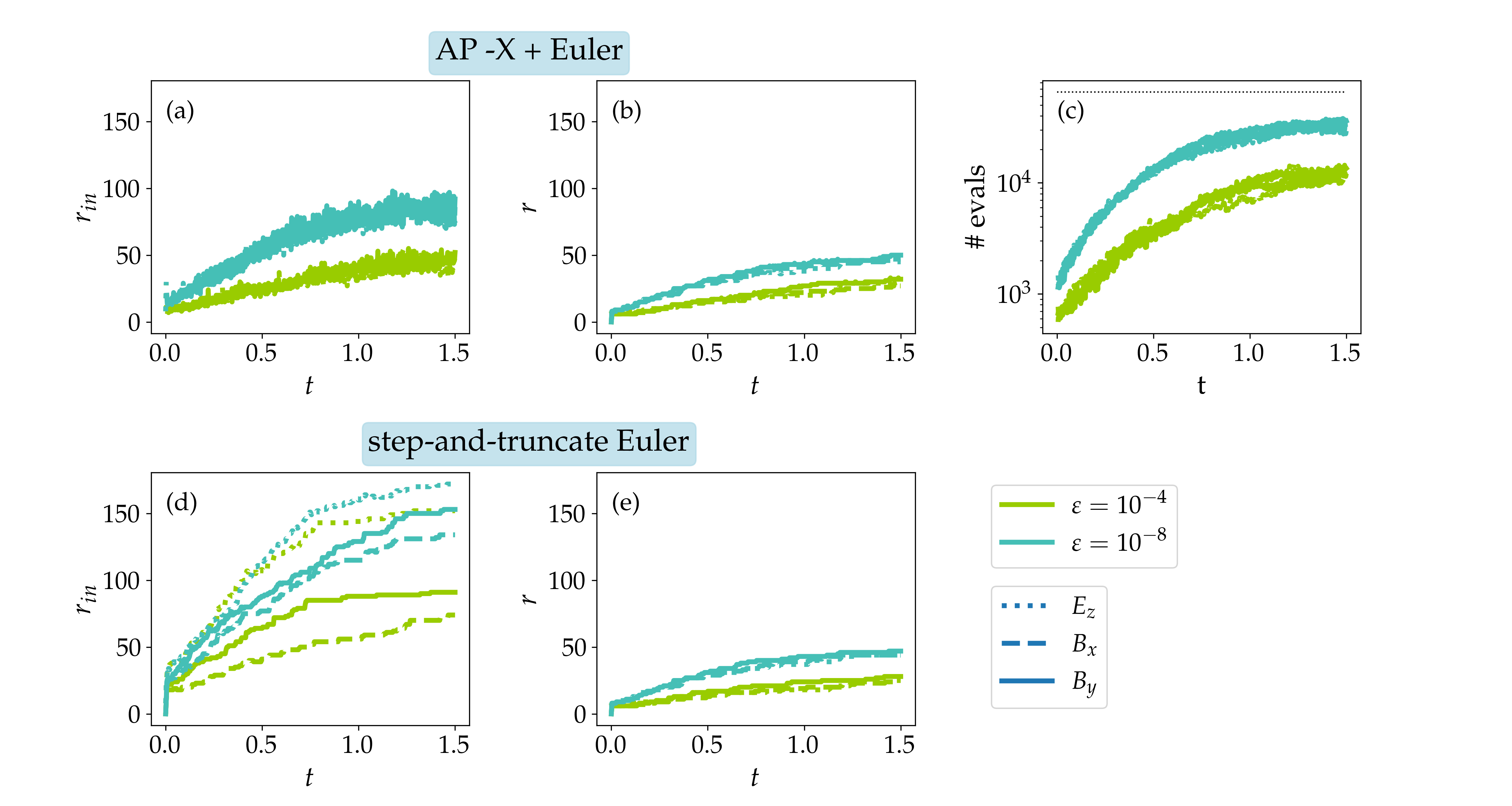}
    \end{subfigure}
    \caption{Results for the Gaussian packet in a dielectric cavity with $2^8$ grid points per dimension. Plots in the top row correspond to results obtained with AP DLR-X, while plots in the bottom row correspond to results obtained with step-and-truncate (SAT). (a-b) QTT ranks for DLRA and SAT before the final rank truncation procedure ($r_{in}$). The rank is smaller for DLR-X than for SAT, which is beneficial since the computational cost scales polynomially with $r_{in}$. (c-d) QTT ranks after the final rank truncation procedure. Notice that the rank for the DLR-X calculation is slightly larger than that of the SAT calculation after times $t=1$. This arises from the approximate nature of the tangent space constructed during the DLR procedure, and can be mitigated by using $\varepsilon_{in} < \varepsilon$ during the DLRA procedure. (e) The number of Euler evaluations in DLR-X compared to the dense grid calculation (the dotted line at $256^2$ evaluations per time step). 
    }
    \label{fig:em2_box}
\end{figure}

Consider the domain $x,y \in [-1,1)$. The system is initialized with $E_z = \exp(-(x^2 + y^2)/100)$ and $B_x=B_y=0$. The cavity is defined by the index of refraction $n$, which varies in space and is given by 
\begin{align}
    & m(x,y) = \frac{1}{4} \left[\tanh(w(x+0.5)) - \tanh(w(x-0.5) \right] \times 
    \nonumber 
    \\ 
    & \hspace{2.5cm} \left[\tanh(w(y+0.5)) - \tanh(w(y-0.5)) \right] 
    \\
    & n^{-2}(x,y) = n^{-2}_o (1 - m(x,y)) + m(x,y)
\end{align}
where $w$ determines the sharpness of the box and $n_o$ is the index of the refraction of the material surrounding the cavity. In the following example, the parameters are set to $w=10$ and $n_o=10$. 

The Maxwell's equations are solved approximately using  a finite volume scheme with first order upwinding:
\begin{align}
    & \partial_t E_z = \frac{c^2}{n^2} \left( \ddx B_y - \ddy B_x \right) + \frac{c}{2} \left(\Delta x \, \ddx^2 + \Delta y \, \ddy^2 \right) E_z
    \\
    & \partial_t B_x = -\ddy E_z + \frac{c \Delta y}{2} \ddy^2 B_x\\
    & \partial_t B_y = \ddx E_z  + \frac{c \Delta x}{2} \ddx^2 B_y 
\end{align}
where a second-order centered finite difference stencil is used to approximate the first and second derivatives.
For simplicity, this work does not carefully treat the inhomogeneity of the material. 
%
%
Additionally, the dissipative upwinding term was derived assuming vacuum material properties. However, for large $n_o$, the fields decay quickly once outside the cavity, justifying the approximate treatment. 
Additionally, the involutions are not satisfied. Achieving this during the DLR procedure in the QTT framework is an area of future research. However, they can be satisfied using divergence correction methods or optimization techniques~\cite{Gourianov2021, einkemmer_low-rank_2020}.

For this example, the AP DLR-X procedure is as follows. 
\begin{enumerate}
    \item blocking: canonicalize $E_z$, $B_y$, and $B_z$. Project the derivative operator onto their respective manifolds. For example, the derivative operator acting on $B_y$ in the time derivative of $E_z$,
    \begin{align}
       {(\ddx)}_\text{eff}^{(B_y \rightarrow E_z)} = \left( \ddx | \textbf{L}_{\a_{i-1}}^{B_y, (i)}, \phi_{\s_i}^{B_y}, \textbf{R}_{\a_i}^{B_y, (i)} \rangle \right)\left[ \mathcal{I}_{\a'_{i-1}}^{E_z, (i)}, \s'_i, J_{\a'_i}^{E_z, (i)} \right] 
       \label{eq:em_ddx_eff}
    \end{align}
    Note that if the field components were in canonical form, \[(\partial_x)_\text{eff} B_{y,\text{eff}} = (\partial_x B_y)_\text{eff} = (\partial_x B_y)[\mathcal{I}^{E_z,(i)},\s_i,\mathcal{J}^{E_z,(i)}] \]
    Similar expressions can be obtained for the remaining operators.
    
    The index of refraction is projected onto the same set of sub-selected indices as $E_z$:
    \begin{align}
       n_\text{eff} = n \left[ \mathcal{I}_{\a'_{i-1}}^{E_z, (i)}, \s'_i, J_{\a'_i}^{E_z, (i)} \right] 
       \label{eq:em_n_eff}
    \end{align}
    
    \item solving: Solve Maxwell's equations at the selected points, using element-wise addition and multiplication.  e.g.,
    \[ E_{z,\text{eff}} (t+\dt) = E_{z,\text{eff}} + \frac{1}{n^2_\text{eff}} \left( (\partial_x)_\text{eff} B_{x,\text{eff}}  - (\partial_y)_\text{eff} B_{y,\text{eff}} \right) \]

    \item decimation: perform the decimation procedure via projection for all field components.
    
\end{enumerate}

Fig.~\ref{fig:em2_box} compares results obtained using the AP method with Euler time stepping to those obtained using a global-update step-and-truncate (SAT) step to compute the Euler time step. The AP method and SAT method yield results of comparable final rank $r$. However, the intermediate $r_{in}$ appears to be slightly larger for the SAT method, corresponding to higher computational costs.

\begin{figure}
    \centering
    \includegraphics[width=0.27\linewidth, trim={0cm 0.0cm 0cm 0cm}, clip]{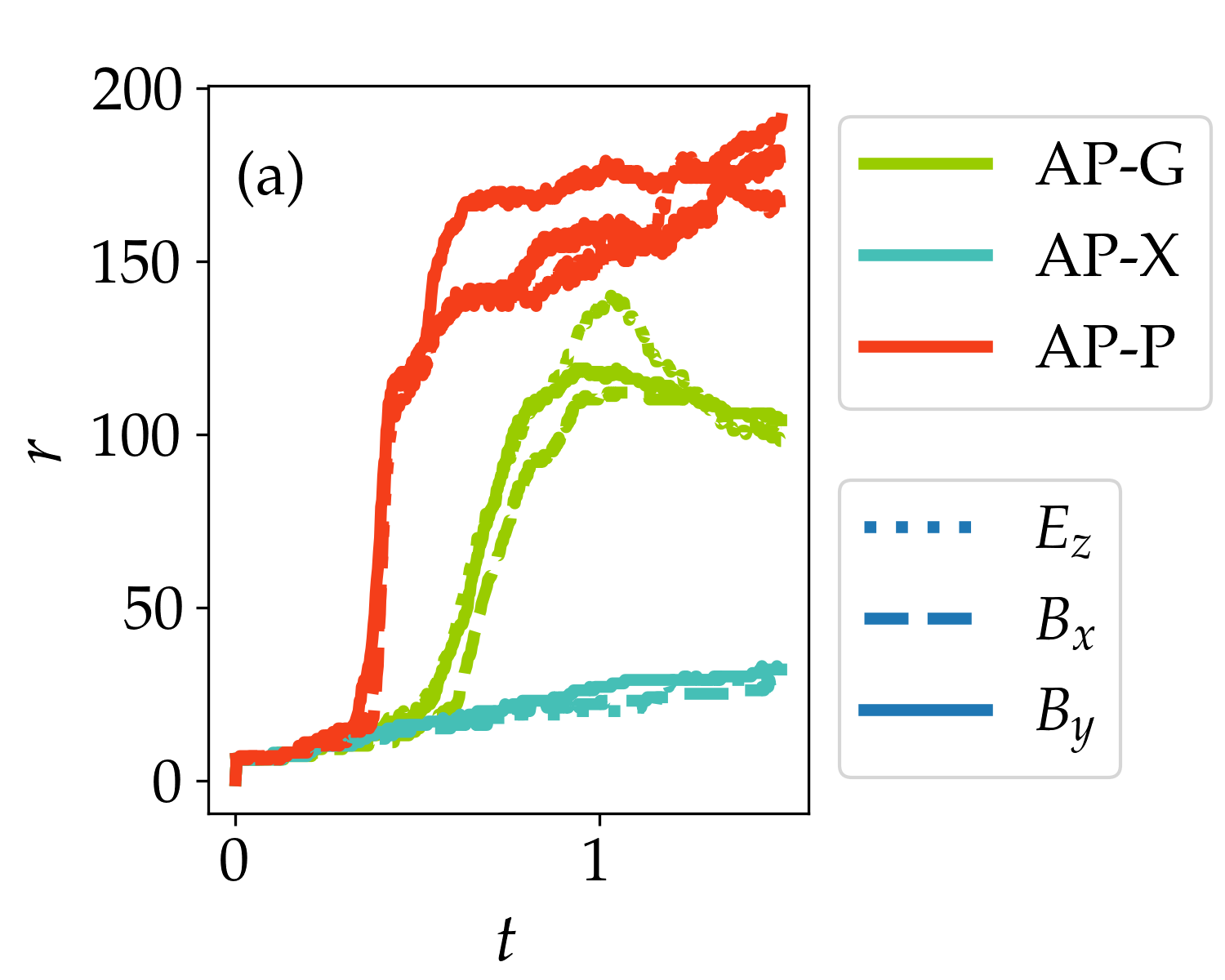}
    \includegraphics[width=0.72\linewidth, trim={0cm 0.2cm 0cm 0cm}, clip]{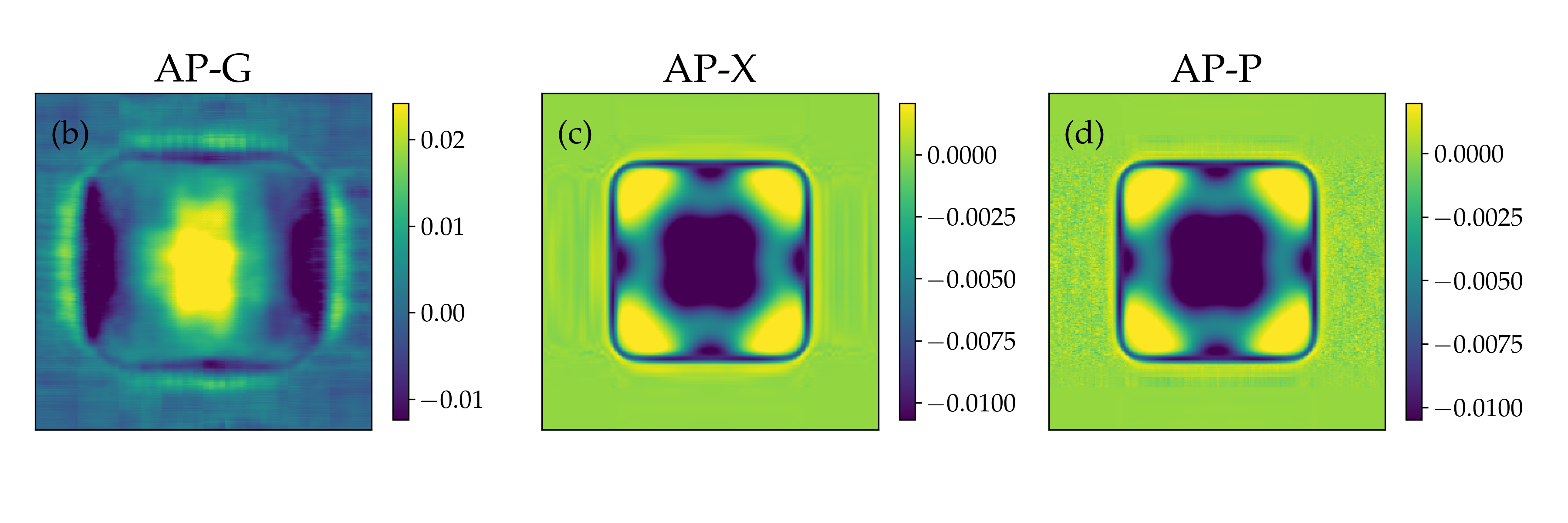}
    \caption{Comparison of results for different AP DLR methods: the standard AP DLR time integration (AP-G), DLR time integration on an interpolative basis (AP-X), and DLR time integration with oblique projection (AP-P). Calculations are performed with $2^8$ grid points per dimension and an error threshold $\varepsilon=10^{-4}$. (a) Ranks of the QTT over simulation time. (b-d) Image of the $E_z$ field at time $t=1.2$ obtained using the three methods. The AP-G result is far from the true solution. The AP-P result shows more noise than the AP-X result.}
    \label{fig:em2_box_compare}
\end{figure}

Despite the element-wise multiplication, a DLR-G procedure can be implemented by treating the element-wise multiplication as the multiplication by a diagonal operator and building an effective operator 
\[ (n^{-2})_\text{eff} = \langle  \textbf{L}_{\a'_{i-1}}^{E_z,(i)}, \phi_{\s'_i}^{E_z}, \textbf{R}_{\a'_i}^{E_z,(i)} \, |\,  (\text{diag}(n^{-2})) \, | \, \textbf{L}_{\a_{i-1}}^{B_y,(i)}, \phi_{\s_i}^{B_y}, \textbf{R}_{\a_i}^{B_y,(i)} \rangle. \]
However, DLR-X requires significantly lower rank than DLR-G and yields more stable results (see Fig.~\ref{fig:em2_box_compare}).

Two DLR-P variants were tested by considering two blocking procedures. The first was to utilize the projectors introduced in Sec.~\ref{sec:oblique}, yielding 
\begin{align}
    & n_\text{eff} = (\textbf{E}_Q^\dagger \textbf{E}_Q \, n)\, [\mathcal{I}^{E_z,(i)}, : ,\mathcal{J}^{E_z,(i)}]
    \\
    & (\partial_x)_\text{eff} = (\textbf{E}_Q^\dagger \textbf{E}_Q \, \partial_x \textbf{E}_Q^\dagger) \, [\mathcal{I}^{E_z,(i)}, : ,\mathcal{J}^{E_z,(i)}]
\end{align}
where $ \textbf{E}_Q $ projects information in the full computational space onto the orthonormal TT manifold, $ | \, \textbf{L}_{\a_{i-1}}^{E_z,(i)}, \phi_{\s_i}^{E_z}, \textbf{R}_{\a_i}^{E_z,(i)} \rangle $.  However, the performance of this method is as bad if not worse than orthonormal DLR. 

The second blocking option is to use the same projection as in DLR-X (Eqs.~\eqref{eq:em_ddx_eff} and ~\eqref{eq:em_n_eff}). This avoids the issue of the accuracy of the mask after projection, allowing for much more stable calculations. However, a large increase in rank is still observed. This is likely because the orthonormal basis captures more of the noisy oscillations, leading to higher ranks as the noise accumulates. 


\subsection{Linear 2-D advection with Boltzmann's equation}

In the final numerical example, consider the behavior of an electron distribution function $f(v_x,v_y)$ in the presence of a magnetic field $B_0 \, \hat{z}$ and time-varying electric field $E_x \hat{x}$. The equation of motion is
\begin{align}
    \partial_t f + \frac{q}{m} \left( ( E_x(t) + B_0 v_y ) \partial_{v_x} - B_0 v_x  \partial_{v_y} \right) f = 0.
\end{align}
With the electric field of the form $E_x(t) = E_0 \cos(\omega t)$ and a Gaussian initial condition, $f(v_x,v_y)= \exp(-(v_x^2 +v_y^2)/2v_{th}^2)$, the resulting dynamics is the drifting of the Gaussian centered at
\begin{align}
    u_x(t) &= \frac{E_0}{B_0} \frac{\omega_c^2}{\omega_c^2 -\omega^2} \left(\sin(\omega_c t) - \frac{\omega}{\omega_c} \sin(\omega t) \right), \\
    u_y(t) &= \frac{E_0}{B_0} \frac{\omega_c^2}{\omega_c^2-\omega^2} \left(\cos(\omega_c t) - \cos(\omega t) \right),
\end{align}
where $\omega_c = q B_0 /m$ and $B_0>0$, and $\omega \neq \omega_c$. In the following calculations, the units are defined such that $B_0=1$, $q/m=-1$, and $v_{th}^2=k_BT/m=1$. The problem parameters are chosen to be $E_0=0.9$, and $\omega=0.4567$.

To accurately measure the time-step error for higher order time integration methods, the equation is solved using a Fourier representation to avoid any grid discretization errors (see Appendix \ref{sec:k} for further details). This choice in basis is less than ideal, since the problem does not have periodic boundary conditions. However, using Hermite polynomials is difficult due to shifting center of the Gaussian, and numerically computing the fast Fourier transform of $v_x$ and $v_y$ on the bounded domain appears sufficient for this problem. 
The Fourier space is discretized with $2^L$ modes along each dimension, leading to a QTT of length $2L$. {The range of frequencies along each dimension is $[-\frac{2^L \pi}{24 v_{th}} , \frac{2^L \pi}{24 v_{th}})$, which corresponds to a real-space velocity grid with the domain $[-12 \, v_{th}, 12 \, v_{th})$.}

\begin{figure}
    \includegraphics[width=\linewidth]{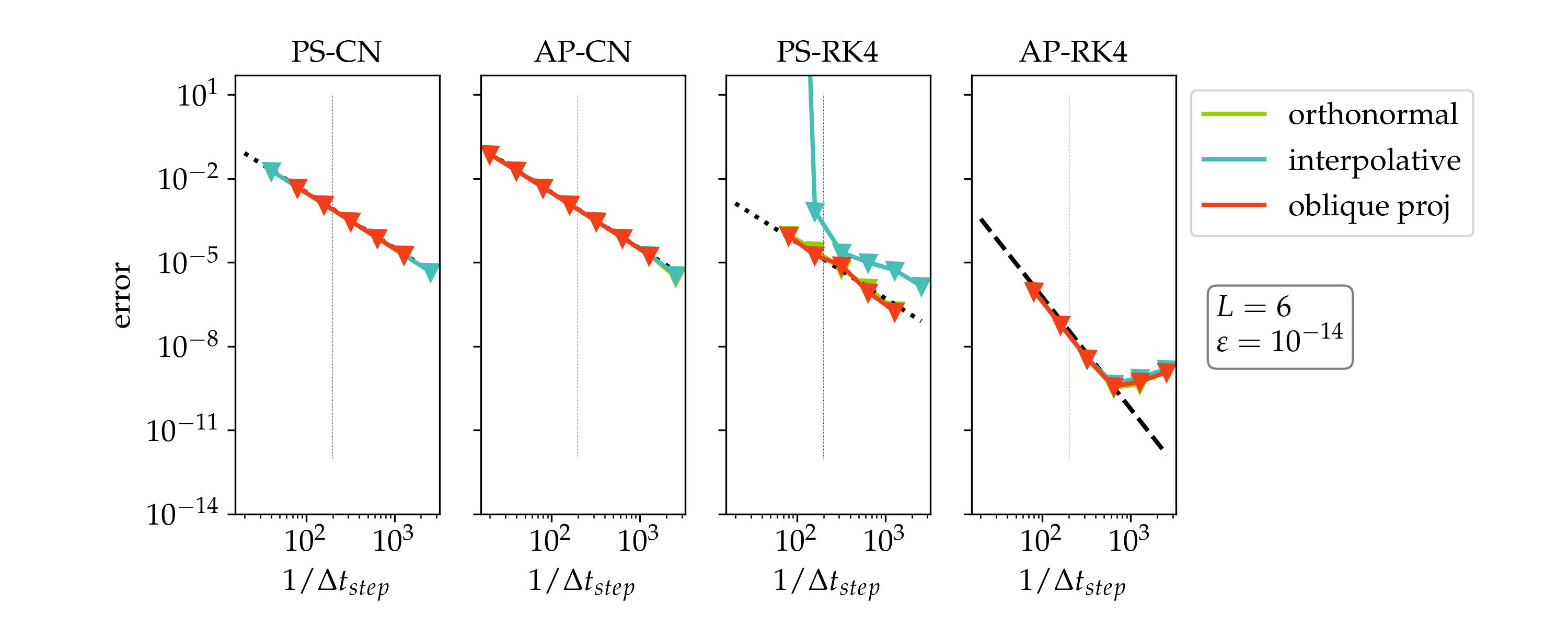}
    \caption{Error at the last time step ($T=100$) with respect to time step size for DLR-type solvers (a) projector splitting with the interpolative construction (PS-X) with RK4 for the internal time integrator (b) PS-X with Crank-Nicolson, (c) basis-update with interpolative construction (BU-X) with RK4 and (d) BU-X with Crank-Nicolson. Different colors correspond to different error measurements: errors in the variance (green) and drift velocity (blue) of the fitted Maxwellian, and the $\ell_2$ error (red) with respect to the theoretical distribution function. The black dashed line denotes $\mathcal{O}(\Delta t^{-4})$ scaling, and the black dotted line denotes $\mathcal{O}(\Delta t^{-2})$ scaling. The thin vertical line denotes the theoretical CFL time step size.}
    \label{fig:advec_exact}
\end{figure}

Because this equation is linear, one can perform time integration using both the orthonormal \cite{ye_quantized_2023} and the interpolative QTT construction. 
The procedure at each tensor core in the DLR-X algorithm looks like
\begin{enumerate}
    \item blocking: determine the manifold from frozen tensor cores of $f$. Project the operator onto the low-rank manifold, e.g.,
    \begin{align}
        (A^{(i)}_{\text{eff},1})_{\a'_{i-1},\s'_i,\a'_i;\a_{i-1},\s_i,
        \a_i} = B_0 \, (v_y \partial_{v_x} |\lenv{i}, \phi_{\s_i}, \renv{i}\rangle )[\mathcal{I}^{(i)}_{\a'_{i-1}}, \s'_i, J^{(i)}_{\a'_i}]
    \end{align}
    \item solving: Compute $M^{(i)} = f^n_\mathcal{I}$ and $A_{\text{eff},1} M^{(i)} = (v_y\partial_{v_x} f^n)_{\mathcal{I}}$ and other remaining terms, yielding the total projected time derivative, $A^{(i)}_\text{eff} M^{(i)}$. Perform time integration. 
    \begin{itemize}
        \item Explicit Euler: Compute $f^{n+1}_\mathcal{I} = M^{(i)} + A^{(i)}_\text{eff} M^{(i)}$.
        \item Multi-stage methods follow directly. However, the tensor $M^{(i)}$ is replaced by the input of each stage. To accurately capture higher order schemes, one must appropriately adjust $A_\text{eff}^{(i)}$ due to the time-dependence of the electric field.
        \item Crank-Nicolson: one needs to solve the linear system \[(\mathbb{I} - \frac{\dt}{2} A^{(i)}_\text{eff} (t_{n+1/2})) f^{n+1}_{\mathcal{I}} = (\mathbb{I} + \frac{\dt}{2} A^{(i)}_\text{eff}(t_{n+1/2})) f_\mathcal{I}^n,\] where $\mathbb{I}$ is the identity. This equation is solved using the conjugate gradient squared method to an accuracy of $10^{-10}$. Note that to achieve second-order accuracy with respect to time, the electric and magnetic fields are evaluated at time $t_{n+1/2} = t_n + \dt/2$.
    \end{itemize}
    \item decimation: shift the center of orthogonality to the next tensor. 
    \begin{itemize}
        \item For the AP method, expand the subspace of $M^{(i)}$ by targeting all accessible intermediate stages as described in Sec.~\ref{sec:final}.  
        \item For the PS method, the procedure is the same as before. 
    \end{itemize}
\end{enumerate}
The dynamics are computed using DLRA with the orthonormal basis (DLR-G), the interpolative basis (DLR-X), and the oblique projector (DLR-P). The $L_2$ error of the last time step with respect to the analytical solution is measured with respect to time step size. Fig.~\ref{fig:advec_exact} shows the results with minimal rank truncation ($\varepsilon=10^{-14}$). In these calculations, all variants of the alternating projection (AP) DLRA algorithm with RK4 yields fourth-order scaling with respect to $\Delta t^{-1}$. AP DLRA with CN yields second order scaling. The projector splitting (PS) DLRA calcucaltions all yield second-order scaling, regardless of the time integrator used to solve the reduced problem. This is because a second-order Lie-Trotter splitting is used.  Interestingly, the error for PS DLR-X with RK4 yields larger errors than its orthonormal counterparts, even though the errors are comparable for the other time integration schemes. 

\begin{figure}
    \includegraphics[width=\linewidth]{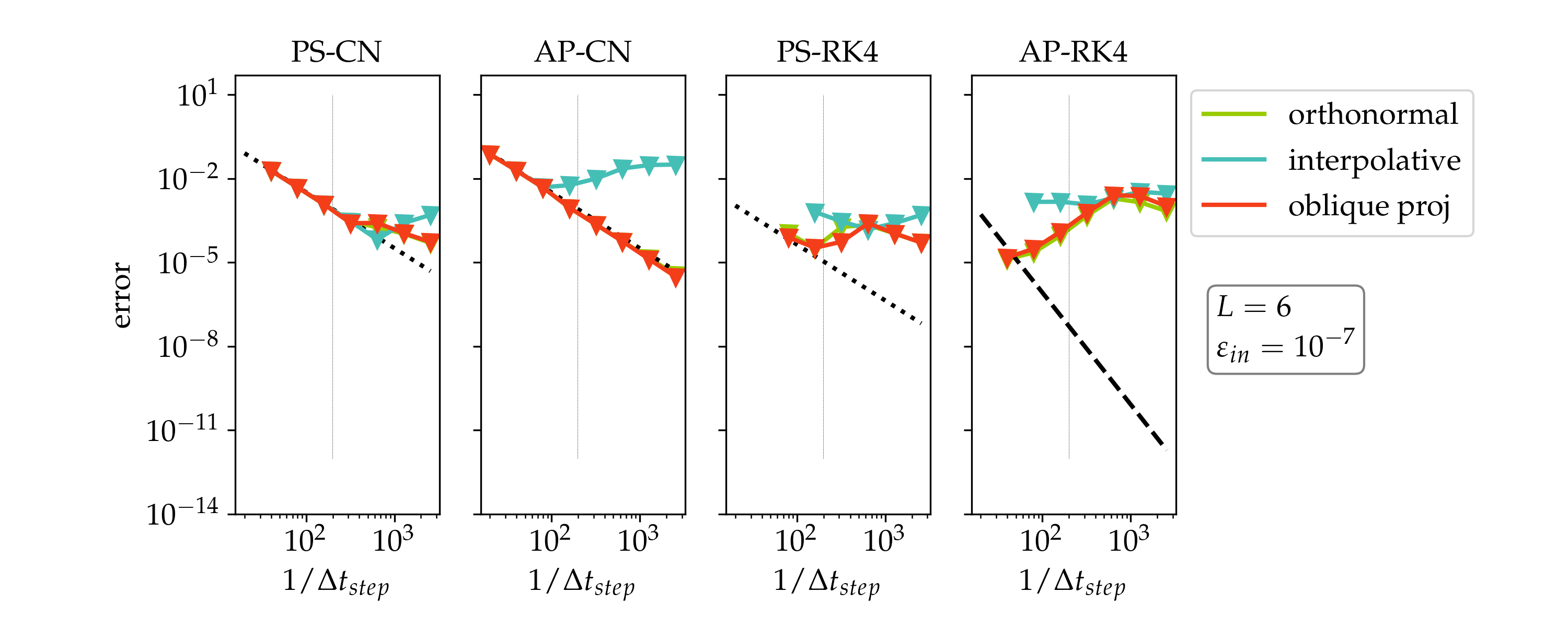}
    \caption{Same as Fig.~\ref{fig:advec_exact} but with an internal cutoff threshold $\varepsilon=10^{-7}$ and final cutoff threshold $\varepsilon_{in}=10^{-6}$. 
    }
    \label{fig:advec_cutoff}
\end{figure}

Results for calculations performed with rank truncation ($\varepsilon=10^{-6}$ and $\varepsilon_{in}=10^{-7}$) are shown in Fig.~\ref{fig:advec_cutoff}.  
With rank truncation, the errors first decay at large time steps but then increases again at smaller time steps. 
Thus, unless higher order time integrators offer improved numerical stability and allow one to take larger time steps, or are required to precisely obtain a low-rank result, their benefits may be lost if the truncation threshold is too large. 

Despite the fact that all of the calculations using the same $\varepsilon$ and $\varepsilon_{in}$, the DLR-X calculations yields larger errors. For example, in the AP-CN calculation, the errors are significantly larger than the errors from DLR-G and DLR-P. The error with respect to $\Delta t$ begins to increase once the error is at about $10^{-2}$, which is much larger the expected error of roughly $\varepsilon$. This is likely because SVD-based rank truncation is not optimal when the tensor train is not in the orthonormal canonical form. 
The performance of the interpolative scheme could potentially be improved if a more optimal rank truncation scheme were used.

\begin{figure}
    \includegraphics[width=\linewidth]{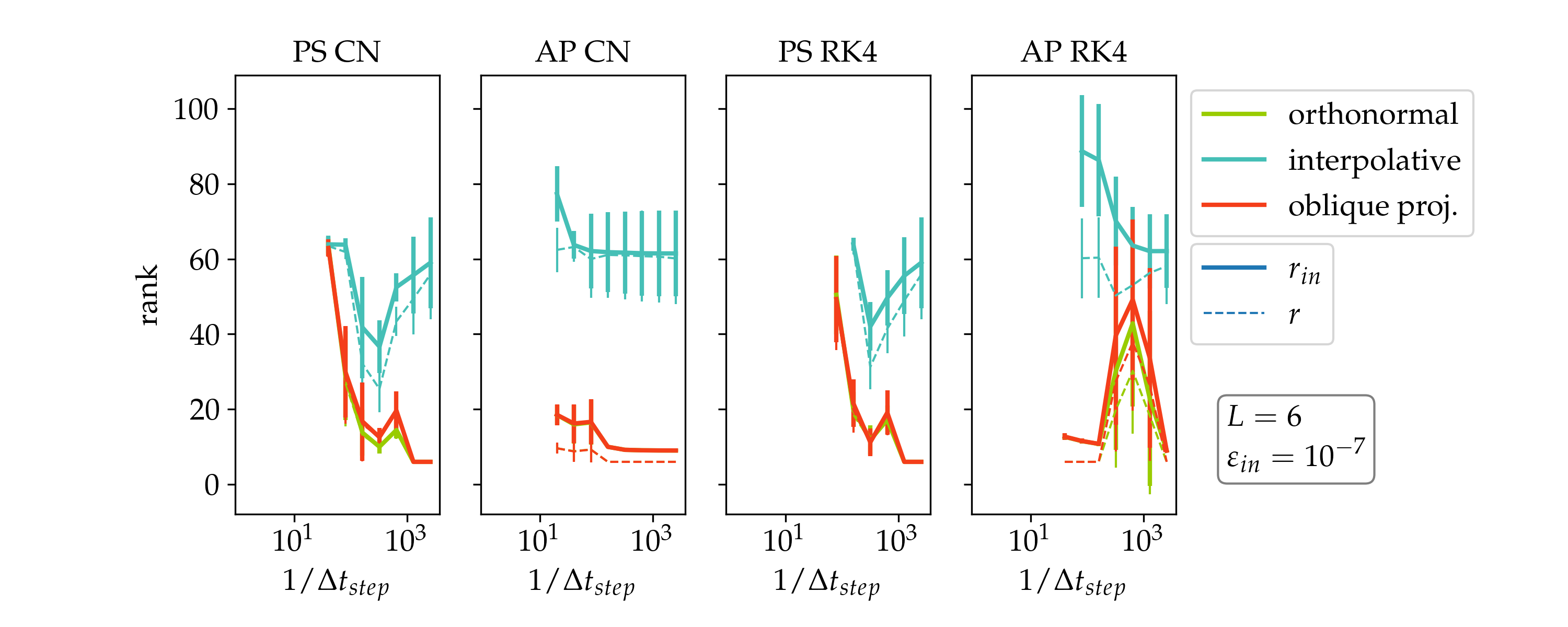}
    \caption{QTT ranks averaged over simulation time for the calculations in Fig.~\ref{fig:advec_cutoff} with $\varepsilon=10^{-6}$ and $\varepsilon_{in}=10^{-7}$. The bars denote the standard deviation of the ranks over the simulation time. Solid lines correspond to the rank during the DLRA procedure ($r_{in}$) and dashed lines correspond to the rank following the rank truncation procedure at each time step ($r$). Ranks for the orthonormal construction (green) are significantly lower than that of the interpolative construction (blue), especially in the case of AP with Crank-Nicolson.}
    \label{fig:advec_cutoff_rank}
\end{figure}

In addition to larger errors, the QTT ranks resulting from DLR-X were also significantly larger, at times appearing to be close to full rank (see Fig.~\ref{fig:advec_cutoff_rank}). The enlarged rank is likely due to the increased numerical noise that occurs in the DLR-X scheme and accumulates over simulation time.


\begin{figure}
    \includegraphics[width=\linewidth]{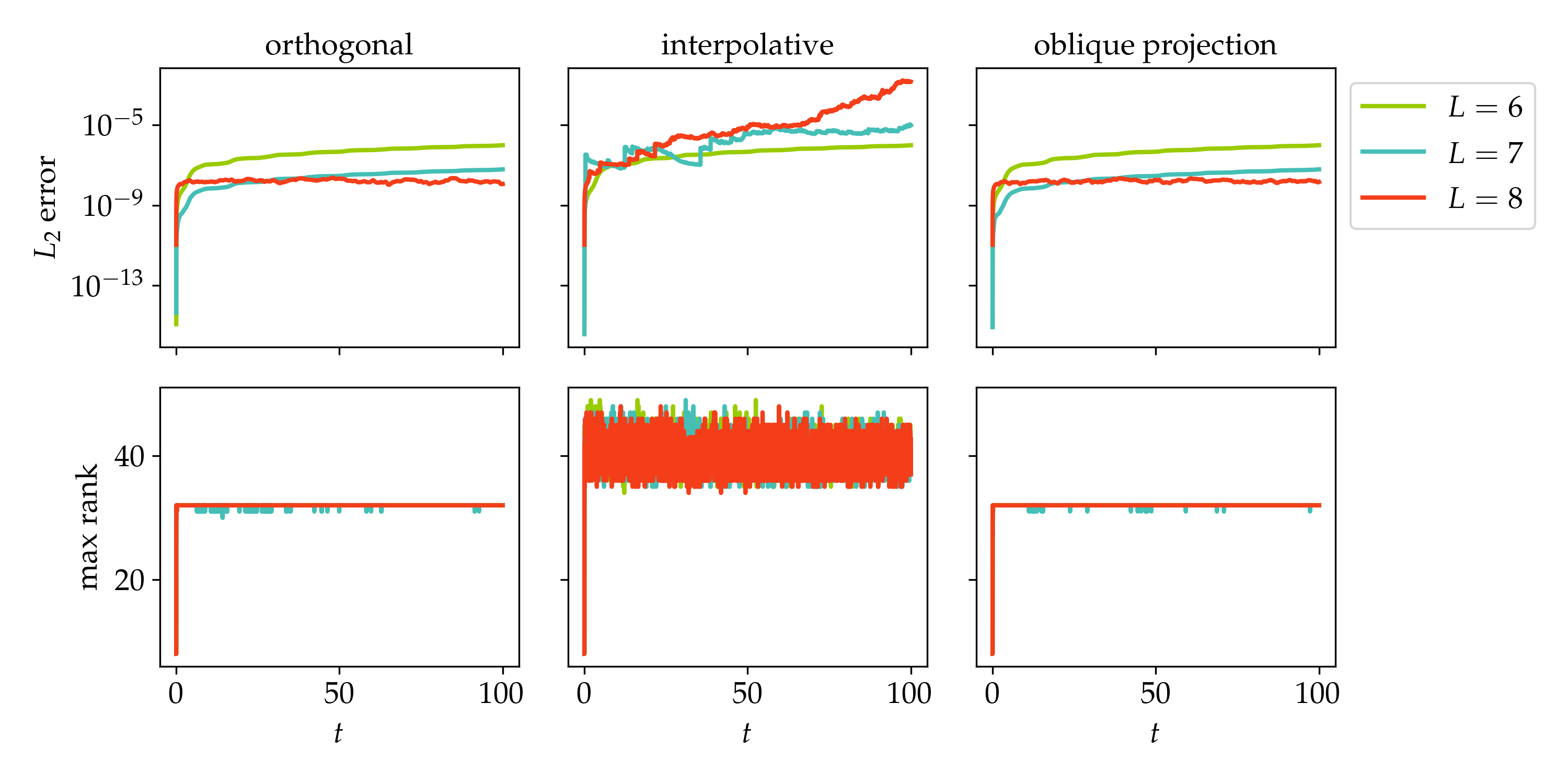}
    \caption{Alternating projection time integration for different resolutions $L$, with fixed rank $r_{max}=8$. The rank requirement is imposed at the end of each time step, and the basis expansion is performed without truncation. (a) $L_2$ errors with respect to the expected theoretical result. (b) QTT rank during the DLR procedure, $r_{in}$. The theoretical maximum is 32, but oversampling in the AP-X procedure improves stability.}
    \label{fig:advec_L}
\end{figure}

Next, consider a fixed rank AP calculation, where the QTT after the DLR time step is compressed to a maximum rank of $r_\text{max}=8$. The error in these calculations are measured for calculations of increasing resolution, and are shown in Fig.~\ref{fig:advec_L}. During the DLR procedure, no rank truncation is performed. For AP-RK4, where the basis expansion targets four stages of the time integration, the maximum rank during the DLR procedure is $4r_\text{max}$. While this does not translate into significant rank reduction for the $L=6$ calculation, the DLR computational costs only increase linearly with $L$, as opposed to exponentially in the dense calculation. For the $L=6$ calculation, a time-step size of $\dt=0.0125$ is used. The time step size is reduced by a factor of two as $L$ is increased by one.

DLR-G and DLR-P yield good results, with the error even decreasing at increased $L$ even though the rank was not increased. However, as resolution was increased, DLR-X struggled with numerical stability. It was observed that by oversampling the CUR decomposition procedure, one obtained better stability. 
The oversampling was accomplished by performing row selection on each of the individual contributions $M^{(i),[k]}$ used in the original basis expansion procedure. These individual contributions were then used to supplement the above indices obtained from the expanded matrix ($I^{all}$), such that the final set of selected rows is their union, $I = I^{all} \cup \left(\bigcup_{k=1}^m I^{[k]} \right) $. In the worst case, $2 m r_i$ indices are selected. However, in practice, it was observed that $I^{all}$ is only expanded by a few indices, still yielding roughly $m r_i$ sampled indices.
Unfortunately, oversampling only mitigated the instability and did not fully resolve the issue.

Another caveat with oversampling is that the tensor train is no longer strictly in canonical form, as it does not satisfy Eqs.~\eqref{eq:canon_cx} and ~\eqref{eq:canon_xr}. Instead, it only satisfies the weaker condition on canonical form, (Eqs.~\eqref{eq:canon_cx_weak} and~\eqref{eq:canon_xr_weak}). If the center of orthogonality is updated beyond a scaling by a constant multiplicative factor, then the TT is no longer in canonical form.  In the case of the AP-X integrator, this is not a severe issue, since canonical form is only lost at the update of the last tensor core. However, this is not the case for the PS-X integrator, as oversampling would cause the QTT to lose canonical form during the sweeping procedure since the orthogonality centers are always updated with a new value.

DLR-X calculations for AP-CN was even more unstable. As such, AP-CN, which would have allowed one to keep the same time-step size regardless of resolution $L$, was not investigated. In contrast, the DLR-G and DLR-P calculations did not face these difficulties.

One potential explanation for the extremely poor performance of DLR-X is because a Fourier representation was used. The representation of a Gaussian is complex, and there is a fast oscillation in the phase of the complex values. It has been observed that the QTT cross approximation tends to struggle with highly oscillatory functions \cite{Lindsey2023multiscale}, so this may be a consistent observation.

\section{Discussion}\label{sec12}

This paper considered three dynamical low-rank time integration schemes. First is the standard procedure (DLR-G) in which the dynamics is projected onto a reduced manifold constructed from orthonormal basis. The projection operation is an orthogonal projection. Second is the analogous procedure that is obtained by considering a TT in an interpolative representation constructed using CUR decomposition (DLR-X). Third is the DLR procedure that converts between the orthonormal and interpolative representations using oblique projections as needed by the problem of interest (DLR-P). 

The above numerical examples demonstrated that each have their strengths and weaknesses: DLR-G is well-suited for linear problems, and typically is very effective. DLR-X appeared to outperform DLR-G in the second application of an electromagnetic field in a dielectric cavity. However, in the case of simple linear advection, the interpolative basis appeared to be a much less efficient representation and prone to numerical issues. This was exemplified in the advection test problem with finite rank. Even though no rank truncation was performed during the DLR sweep, the limited rank of the TT manifold even after basis expansion (prior to oversampling) resulted in unstable calculations. DLR-P offers a middle ground. Its flexibility allows one to treat nonlinear element-wise operations (which are common for upwind time integration schemes) while retaining the efficiency of the orthonormal basis. If implemented properly, DLR-P offers performance somewhere between DLR-G and DLR-X.

The computational cost of the each DLRA time step scales like $\mathcal{O}(dr^2_{in} r\log(N))$, where $r$ is the tensor train rank after compression, $r_{in}$ is the tensor train rank during the DLRA procedure, $d$ is the size of the quantized dimensions, and $N$ is the total number of grid points, which scales exponentially with problem dimension. The cubic scaling with respect to the rank (assuming $r_{in}$ and $r$ are comparable) is a significant overhead, and calculations requiring full rank will be much slower than the standard dense calculation. The rank required for a given calculation is highly problem dependent, and will often change over time for advection-based problems.  In general, it is difficult to predict the rank required to achieve the desired numerical accuracy. This work does not make any attempt to demonstrate the efficiency of QTTs over a dense vector representation. This is in part because these test problems are quite small. One could choose to increase the resolution of the problem to demonstrate the utility of QTTs, as was done in the Burger's equation example. However, in practice, if one could avoid increasing resolution, one would. As such, it is expected that QTTs will exhibit most advantage in multi-scale problems requiring high resolution. QTTs can also reduce the computational cost of multi-dimensional problems, though non-quantized tensor trains will often suffice unless high resolutions are needed.  
Clearly demonstrating the advantage of QTTs is an active area of investigation.

Though not addressed in this work, 
one challenge low-rank methods face is the loss of conservation of certain quantities of interest. Several methods have been discussed in the context of non-quantized tensor trains for specific problems (e.g., the Vlasov equation \cite{einkemmer_robust_2023, coughlin_robust_2024, guo_conservative_2024}). One solution is to consider a micro-macro decomposition, such that moments of particular interest (such as mass, momentum, and energy) are treated exactly, while others are treated approximately with this low-rank method. 
However, if one insists on a conservative QTT method, one would need to ensure conservation in two places. First, in the DLRA time step, one must ensure that the updated TT manifold can accurately capture the conserved quantities. Second, in the rank truncation step, one must ensure that the no truncation occurs on any part that will affect the conserved quantities. For some spectral methods, in which the conserved quantity is solely the coefficient of a specific spectral mode, conservation may be more naturally achieved in the interpolative construction. 
In other cases, such as in real-space and Fourier-space, an orthonormal basis construction is more amenable to enforcing conservation. Having the option to convert between the interpolative and orthonormal representations will be helpful in these scenarios.

\section{Conclusion}

The interpolative TT construction allows one to interpret the tensor train as a polynomial interpolation with degree $r$ polynomials.
This work demonstrates that interpolative DLR-type time integration can be used successfully in the context of quantized tensor trains, allowing one to solve problems requiring element-wise operations that may appear as a result of choosing an upwind time integration scheme, or as a result of the equations of motion themselves. With a well-chosen basis expansion procedure, one can obtain higher-order time integration schemes, though the benefits higher-order methods are negated when using a modest truncation threshold since the low-rank approximation error will dominate. 

The success of the DLRA time integration scheme is highly problem-dependent. While the interpolative TT construction yields better results for some types of problems, it can also lead to numerical instabilities in other scenarios. This issue can be mitigated by using an oblique projector to project the interpolative construction onto a appropriately chosen orthonormal construction. In doing so, one can relate interpolative DLR to the standard DLR methods, which tend to be stable but limited to linear problems.
This paper provides a practical survey of interpolative DLRA time integration schemes and demonstrates how one can now consider a wider range of dynamical systems.

\backmatter





\bmhead{Acknowledgements}
This material is based in part upon work supported by the Laboratory Directed Research and Development Program of Lawrence Berkeley National Laboratory under U.S. Department of Energy Contract No. DE-AC02-05CH11231 (E.Y.), and in part upon work supported by the U.S. Department of Energy, Office of Science, Office of Advanced Scientific Computing Research, Scientific Discovery through Advanced Computing (SciDAC) program through the FASTMath Institute. This research used resources of the National Energy Research Scientific Computing Center, a DOE Office of Science User Facility supported by the Office of Science of the U.S. Department of Energy under Contract No. DE-AC02-05CH11231 using NERSC award ASCR-ERCAPm1027.
E. Ye would like to acknowledge A. Dektor for helpful discussions regarding his work on interpolative DLR methods.

\section*{Statements and Declarations}


\begin{itemize}

\item Competing interests. 

The authors declare no competing interests.

\item Code availability 

All code was implemented in Python using the tensor network package quimb \cite{gray2018quimb}. The source code for the above method will be made available on Github.

\end{itemize}

\noindent

\begin{appendices}

\section{QTT in Fourier space} \label{sec:k}
In Section 2.3, the PDE is solved in the Fourier representation instead of the standard real-space basis. This is still efficient in the low-rank representation. 
In Fourier space, the quantization procedure is the same, except the basis functions are now of the form $ \phi_n(x) = \exp(i2\pi n x) $, and a binary mapping between $n$ and hyperindices ($\s_1,\hdots \s_L$) are used:
\[n = \sum_{k=1}^L 2^k \sigma_k = \sigma_L \hdots \sigma_3 \sigma_2 \sigma_1.\]
This choice of basis still offers a multiscale representation, with the splitting of the $\sigma$ at grid level $m$ separating the low frequency modes $(\sigma_1,...\s_{m})$ from the high frequency modes $(\sigma_{m+1}, ... \s_L)$
\[
n_{\leq m} = \sum_{k=1}^m 2^k \sigma_k, \quad n_{> m} = \sum_{k=m+1}^{L} 2^k \sigma_k = 2^m \sum_{k=1}^{L-m} 2^k \sigma_k
\]

Converting the real-space operations in the PDE to Fourier space is done by taking the Fourier transform. In the case that an analytical form is not readily accessible, one can efficiently compute it numerically in the QTT representation since the Fourier transform is low rank \cite{chen_quantum_2023, chen_direct_2024}. Constructing the operator explicitly can also be preferable to avoid issues with artificial ringing effects.

For example, consider the convolution operation, which takes two inputs and generates an output. Ref.~\cite{kazeev_multilevel_2013} explains the construction of the QTT representation, determining that it is of rank 2.
It can also be expressed explicitly \[ C_{ijk} = \sum_{a,b,c} \mathcal{F}_{kc} \delta_{abc} \mathcal{F}_{ai}^{-1} \mathcal{F}_{bj}^{-1} \]
where $\mathcal{F}$ is the matrix representation of the Fourier transform, such that
\[ (x * y)_k = \sum_{ij} C_{ijk} x_i y_j \]
Since the Fourier transform can be written in QTT form, then this operator can be also be efficiently constructed in a way analogous to the zip-up algorithm \footnote{see \href{https://tensornetwork.org/mps/algorithms/zip_up_mpo/}{tensornetwork.org}}: 
starting from $i=1$, contract all tensors associated with index $\sigma_i$. If not at the end of the chain, perform the QR decomposition, setting $Q$ as the $i^\text{th}$ tensor core and associating $R$ with the $i+1$ site. Once the new QTTs are constructed, perform the standard canonicalization procedure sweeping from site $L$ to site $1$.
With an rank-truncation threshold of $\varepsilon = 10^{-14}$, it also returns a QTT of rank 2, regardless of resolution. 

\begin{figure}
    \includegraphics[width=0.45\linewidth]{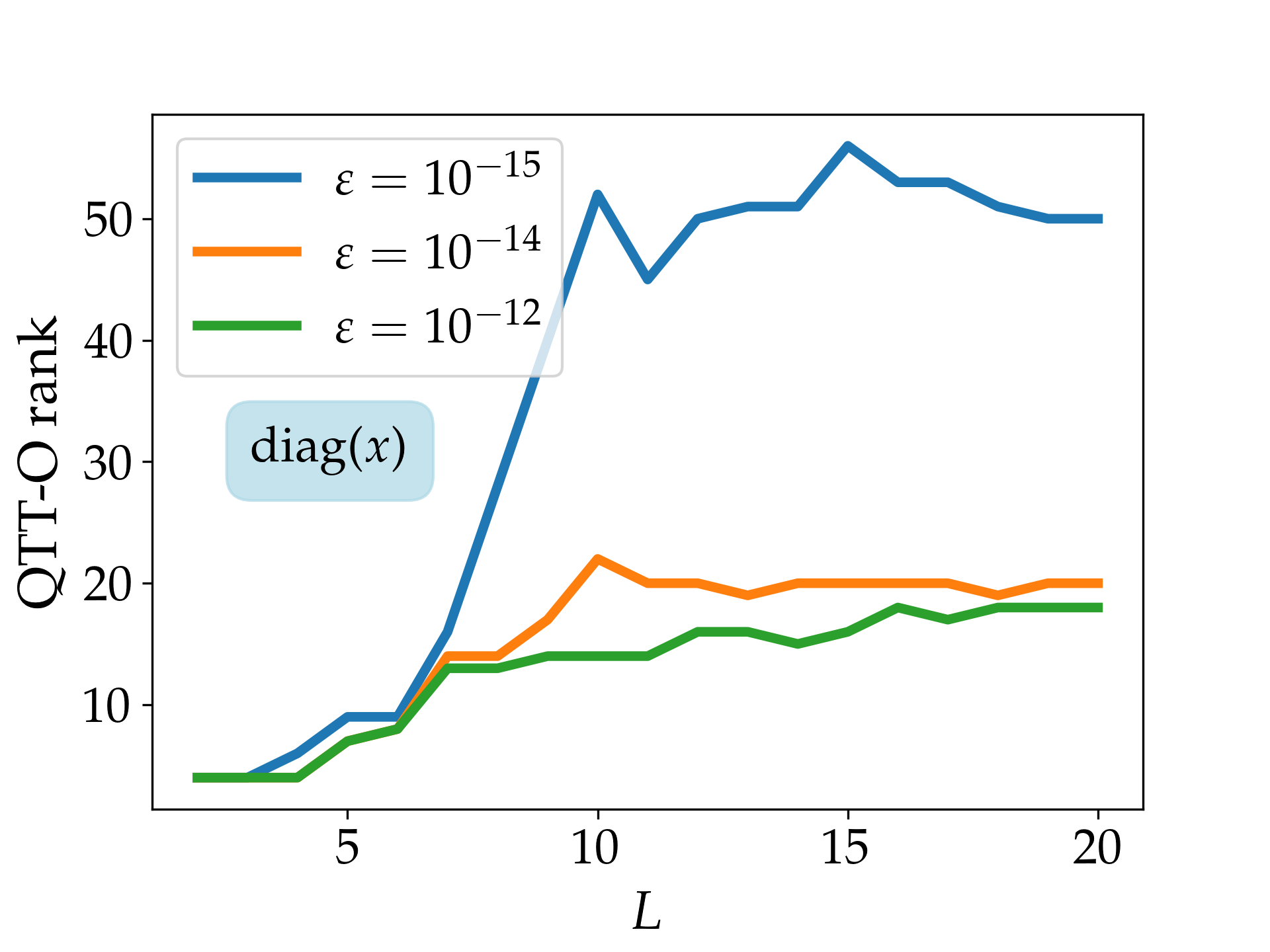}
    \includegraphics[width=0.45\linewidth]{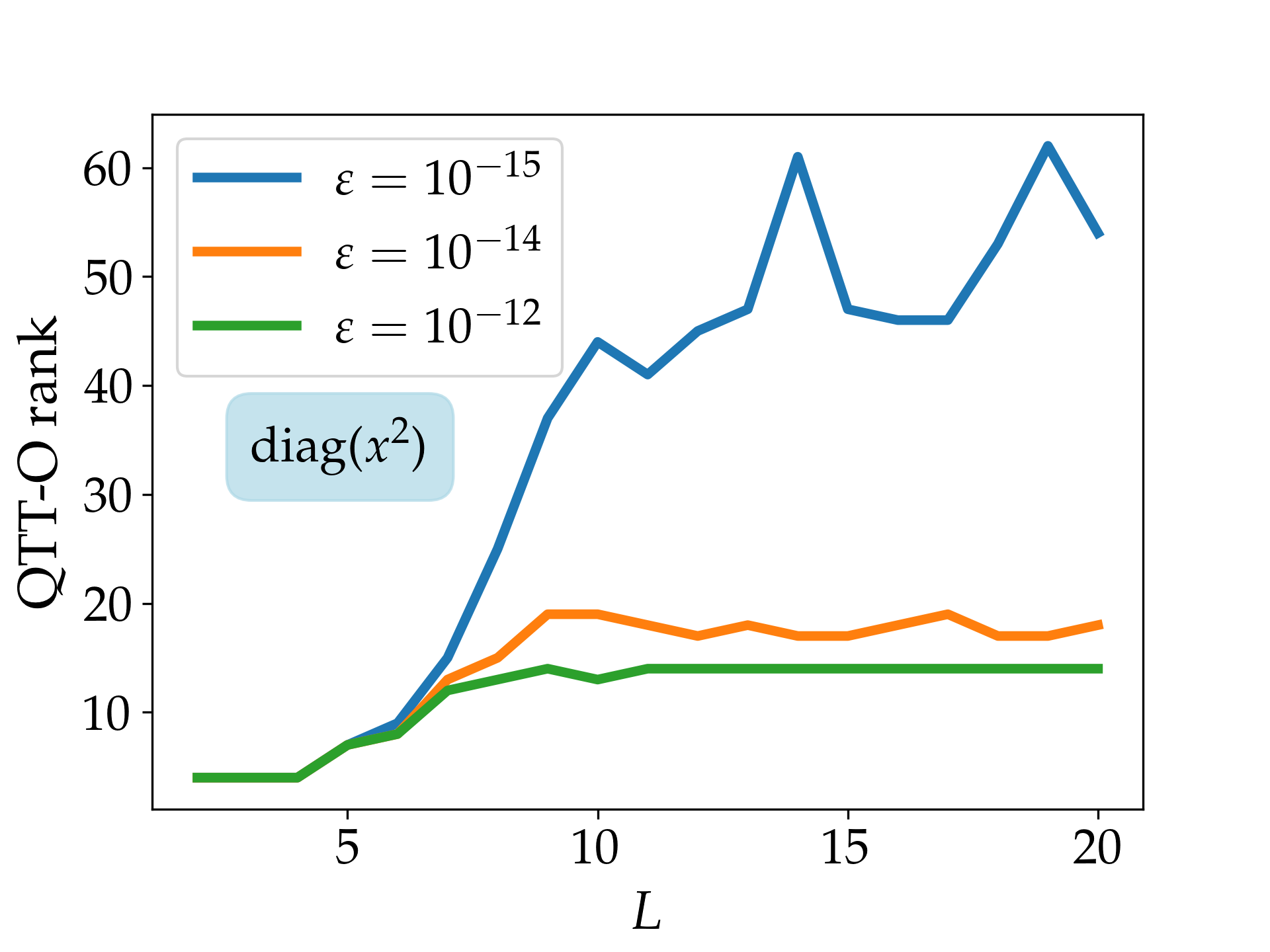}
    \caption{QTT ranks versus $L$ for Fourier representation of $\text{diag}(x)$ and $\text{diag}(x^2)$.}
    \label{fig:k_rank}
\end{figure}

To measure the $m^\text{th}$ spatial moment of some function $\tilde{f}$ that it is in the Fourier basis, one will need to perform the convolution $\mathcal{F}x^m * \tilde{f}$. The linear operator that performs this operation $M \tilde{f} = (\tilde{f} * \mathcal{F}x^m) $ is
\[ M_{ij} = \sum_{a,b,c} \mathcal{F}_{kc} \delta_{abc} \mathcal{F}_{ai}^{-1} \sum_j \left( \mathcal{F}_{bj}^{-1} (\mathcal{F} x^m)_j \right) = \sum_j C_{ijk} (\mathcal{F} x^m)_j \,. \]
For the first and second orders, the ranks are larger than their real-space counterparts, but plateaus at about 20 at high resolutions (large $L$), as shown in Fig.~\ref{fig:k_rank}.

Once these operators are constructed offline, one can solve the PDE in the Fourier basis without explicitly converting to a real-space representation.


\end{appendices}


\bibliography{sn-bibliography}

\end{document}